\documentclass[a4paper,12pt,english]{article}

\usepackage{lineno,hyperref}
\usepackage[english]{babel}
\usepackage{graphicx}
\usepackage{amsmath}
\usepackage{amsfonts}
\usepackage{amssymb}
\usepackage{color}
\usepackage{algorithm}
\usepackage{algorithmic}
\usepackage{authblk}
\usepackage{fancyhdr}
\usepackage{booktabs}
\usepackage{multirow}
\usepackage{amsmath}
\usepackage{url}
\modulolinenumbers[5]

\numberwithin{equation}{section}
\newcommand{\bm}[1]{\mathbf{#1}}
\newcommand{\e}{{\mathrm{e}}}
\newcommand{\alphab}{\boldsymbol{\alpha}}
\newcommand{\betab}{\boldsymbol{\beta}}
\newcommand{\R}{{\mathbb R}}
\newcommand{\opcoef}{\texttt{opcoef}}
\newcommand{\opevmat}{\texttt{opevmat}}
\newcommand{\legmoms}{\texttt{legmoms}}
\newcommand{\chebtwomoms}{\texttt{cheb2moms}}
\newcommand{\dchebmoms}{\texttt{dchebmoms}}

\providecommand{\keywords}[1]{\textbf{\textit{Keywords-}} #1}

\title{Fast and accurate computation of orthogonal moments for texture analysis}
\author[1]{Cecilia Di Ruberto}
\author[2]{Lorenzo Putzu} 
\author[1]{Giuseppe Rodriguez}
\affil[1]{Department of Mathematics and Computer Science, 
University of Cagliari, via Ospedale 72, 09124 Cagliari, Italy}
\affil[2]{Department of Electrical and Electronic
Engineering, University of Cagliari, piazza d'Armi, 09123 Cagliari, Italy}

\begin{document}

\maketitle

\begin{abstract}
In this work we describe a fast and stable algorithm for the computation of the
orthogonal moments of an image.
Indeed, orthogonal moments are characterized by a high discriminative power,
but some of their possible formulations are characterized by a large
computational complexity, which limits their real-time application.
This paper describes in detail an approach based on recurrence relations, and
proposes an optimized Matlab implementation of the corresponding computational
procedure, aiming to solve the above limitations and put at the community's
disposal an efficient and easy to use software.
In our experiments we evaluate the effectiveness of the recurrence
formulation, as well as its performance for the reconstruction task, in
comparison to the closed form representation, often used in the literature.
The results show a sensible reduction in the computational complexity,
together with a greater accuracy in reconstruction.
In order to assess and compare the accuracy of the computed moments
in texture analysis, we perform classification experiments on six well-known
databases of texture images. Again, the recurrence formulation performs better
in classification than the closed form representation. More importantly, if
computed from the GLCM of the image using the proposed stable procedure, the
orthogonal moments outperform in some situations some of the most diffused
state-of-the-art descriptors for texture classification.
\end{abstract}

\section{Introduction}\label{sec:intro}

Many texture feature extraction methods have been reported in the literature.
Among them, the statistical methods compute simple statistical properties
(e.g., first order parameters) or sophisticated properties (e.g., co-occurrence
matrices \cite{Gelzinis,Gong,Sarah91}, or gray level run length \cite{Belur91})
of images. Other texture feature extraction approaches are based on models
(Markov random field and fractals, or local binary pattern
\cite{Ojala96}), mathematical morphology \cite{Soille2002}, or classical
transform methods (Fourier, Gabor, wavelet transforms \cite{Mallat2002}).

The application of a texture analysis method consists of extracting a
set of parameters from an image, in order to characterize the texture it
contains \cite{DiRuberto2015a}. Each parameter expresses a special property,
such as coarseness, homogeneity, or the local contrast. The texture feature
extraction can be carried out either at the pixel level, by calculating the
textural parameters on a very small neighborhood, or at the level of a region
of interest defined by the user, which may correspond to a large number of
pixels.

Moments are statistical measures which can be used to obtain relevant
information on an object. Since Hu introduced them in image analysis
\cite{Hu1962}, moment functions have been widely used as discriminative
descriptors in image processing and pattern classification applications. In
1994, Tuceryan discussed texture feature extraction and texture analysis based
on geometric moments \cite{Tuceryan1994}. Begum applied complex moments to
texture segmentation \cite{Bigun1994}. However, both geometric and
complex moments contain redundant information and are sensitive to noise. This
is due to the fact that the associated kernel polynomials are not orthogonal.
Many authors proposed the use of orthogonal moments: Mukundan adopted the
discrete Chebyshev moments \cite{Mukundan2001}, Yap studied another set of
discrete moments, known as Krawtchouk moments \cite{Yap2003}, Teague suggested
the use of Legendre and Zernike moments \cite{Teague1980}.

Orthogonal moments are shown to be less sensitive to noise and have an
efficient capability of feature representation. They allow to reconstruct the
image intensity function analytically, from a finite set of moments, using
the inverse moment transform. Legendre and Zernike moments are most widely used
because of their minimum redundancy.
Indeed, they can represent the properties of an image with no redundancy or
overlap of information between the moments \cite{Hwang2006}. Because of these
important features, Zernike moments have been widely used in different
types of applications \cite{DiRuberto2005,Li2009}. They have been
utilized in shape-based image retrieval \cite{DiRuberto2008}, edge detection
\cite{Li2010}, and as a feature set in pattern recognition \cite{Haddadnia2003}.
In \cite{Hitam2014}, the authors apply Zernike moment features to
retrieve binary and gray level images from MPEG-7 and COIL-20
dataset, respectively, showing their suitability for image retrieval, due to
their rotation invariance.
A mixture of ternary patterns and Zernike moments has been proposed to provide
an efficient image classification which is robust against illumination and rotation \cite{Obulakonda2016}.

Unfortunately, these approaches are often characterized by a large
computational complexity, making them
unsuitable for real-time applications. This lead many researchers to 
develop faster algorithms for computing Zernike moments.
In \cite{Tahmasbi2011}, the authors propose a CAD system for the diagnosis of
breast masses in mammography images which uses Zernike moments for extracting
the shape and margin properties of the masses.
In \cite{Oujaoura2014}, an image annotation system, based on different type of
moments, has been developed to allow searching image databases. 
The experimental results showed that the annotation system coupling Legendre
moments to Bayesian networks gives good results for images that are well and
properly segmented.

In \cite{Wu2010}, a moment-based approach is proposed for texture analysis of
medical images, namely, CT liver scan and prostate ultrasound.
The neighborhood of a texture pixel is calculated by different moments for
texture feature extraction. After being verified on Brodatz textures, the
moment-based texture analysis method is applied to CT liver scan classification
and prostate ultrasound segmentation. A support vector machine and a
multi-channel active contour model are used in this application. The results
show that the proposed method is promising, but still have some limitations.
Local binary pattern and Legendre moments have been used as features for CT
liver images classification \cite{Vijayalakshmi2016}.
In \cite{ma2011}, a new set of rotation and scale invariants of Legendre
moments is introduced, achieving good results in image classification
experiments.
Lakshmi et al.~\cite{Lakshmi} used Legendre moments for palm print
authentication with a very good prediction accuracy.

The computation of Legendre moments is in general a time consuming process.
In many references \cite{Oujaoura2014,Wu2010,Vijayalakshmi2016,ma2011},
their computation has been performed using closed form representations for
orthogonal polynomials, and taking little care to the accuracy of the
quadrature formulas used to approximate integrals.
In this work, we describe a fast and stable algorithm for the computation of
the orthogonal moments of an image with respect to both a continuous and a
discrete inner product, based on classical recurrence relations for orthonormal
polynomials.
Despite their recurrent structure, such relations can be coded without
recursive calls, allowing us to develop an efficient and accurate 
Matlab toolbox.
We propose this software as a standard tool for orthogonal moments computation
and we make it freely available to the scientific community.
The reconstruction performance of the proposed methods and their discriminative
ability in texture classification are investigated in comparison to orthogonal
moments computed from a closed form representation, and to other classification
methods.
We finally discuss the effect of weighted orthogonal moments on the processing
of two particular datasets.

Our numerical experiments confirm the well known fact that, while recurrence
relations for orthogonal polynomials produce reliable results, closed form
representations should be avoided for moments computation. 
Moreover, computing moments from Gray Level Co-occurrence Matrices (GLCM)
produces a more accurate classification and, at least in some situation,
weighted moments may lead to a slight performance improvement. Orthogonal
moments from GLCM appear to be competitive with state-of-the-art approaches
based on Convolutional Neural Networks (CNN). We give a possible interpretation
for the different performance of orthogonal moments and CNN when applied to
datasets with specific features.

The paper is organized as follows. In Section~\ref{sec:legmom} we introduce the
definition and properties of Legendre moments. The algorithms for the
computation of orthogonal moments contained in the Matlab toolbox are described
in Section~\ref{sec:comp}.
Section~\ref{sec:results} presents the results of numerical experiments
assessing the reconstruction performance and the discriminative accuracy
in the texture analysis of six dataset collections.
Finally, in Section~\ref{sec:conclusion} we summarize the content of the paper
and describe our plans for future work.

\section{Orthogonal moments}\label{sec:legmom}

The use of moments for image analysis and pattern recognition was inspired by
Hu \cite{Hu1962}. Among all types of moments, orthogonal moments present
the peculiar property of being characterized both by small information
redundancy and by high discriminative power.
A representative family of orthogonal moments is the well known Legendre
moments. They were first introduced in image analysis by Teague
\cite{Teague1980} and have been exhaustively used in many pattern recognition
and features extraction applications, due to their invariance to scale,
rotation, and reflection change \cite{chong2004,ma2011}.

Here, we briefly recall some definitions and properties of orthogonal
polynomials; see~\cite{gautschi2004}.
Let us consider the space $L_w^2[-1,1]$ endowed with the weighted inner product
\begin{equation}\label{innprod}
\langle f, g \rangle = \int_{-1}^1 f(x) g(x) w(x)\,dx
\end{equation}
and the induced norm $\|f\|=\sqrt{\langle f, f \rangle}$, being $w(x)\geq 0$ a
weight function.
Two functions $f,g\in L_w^2[-1,1]$ are said to be orthogonal if $\langle f, g
\rangle=0$; they are orthonormal if, additionally, both functions have unitary
norm.

A fundamental theorem states that there exists a unique infinite sequence of
orthogonal polynomials, that is, polynomials $p_k(x)$ of degree $k$, with
$k=0,1,\ldots$, such that $\langle p_k, p_\ell \rangle=0$ whenever $k\neq\ell$.
Orthogonal polynomials may be scaled in different ways, i.e., so that they
are monic (the monomial of maximum degree has unitary coefficient) or so that
they are orthonormal. We will adopt the latter normalization.

Monic orthogonal polynomials associated to inner product \eqref{innprod}
are defined by the following three-term recurrence relation
\begin{equation}\label{monicrecurrence}
\begin{cases}
\tilde{p}_{k+1}(x) = (x-\alpha_k)\tilde{p}_k(x) - \beta_k\tilde{p}_{k-1}(x), \\
\tilde{p}_0(x)=1, \quad \tilde{p}_{-1}(x)=0, 
\end{cases}
\end{equation}
for $k=0,1,\ldots$, where $\tilde{p}_{-1}(x)$ serves the only purpose of
starting the recursion,
$$
\alpha_k=\frac{\langle x\tilde{p}_k,\tilde{p}_k \rangle}{\langle
\tilde{p}_k,\tilde{p}_k \rangle},
\ k=0,1,\ldots, \qquad
\beta_k=\frac{\langle \tilde{p}_k,\tilde{p}_k \rangle}{\langle
\tilde{p}_{k-1},\tilde{p}_{k-1} \rangle},
\ k=1,2,\ldots, \quad
$$
and $\beta_0=\langle \tilde{p}_0,\tilde{p}_0 \rangle$.
This process is often referred to as the Stieltjes
procedure~\cite{gautschi2004}.
Polynomials which are orthonormal with respect to \eqref{innprod}, that is
$p_k(x)=\tilde{p}_k(x)/\|\tilde{p}_k(x)\|$ for $k=0,1,\ldots$,
satisfy the three-term recurrence 
\begin{equation}\label{recurrence}
\begin{cases}
\sqrt{\beta_{k+1}}p_{k+1}(x) = (x-\alpha_k)p_k(x) - \sqrt{\beta_k}p_{k-1}(x), \\
p_0(x)=\sqrt{\beta_0^{-1}}, \quad p_{-1}(x)=0, 
\end{cases}
\end{equation}
with the same coefficients $\alpha_k$ and $\beta_k$ used in
\eqref{monicrecurrence}.

Legendre polynomials correspond to the weight function $w(x)=1$. In this case,
$$
\alpha_k=0, \ k=0,1,\ldots, \quad
\beta_k=\frac{k^2}{4k^2-1}, \ k=1,2,\ldots, \quad
\beta_0=2.
$$
Second kind Chebyshev polynomials are obtained by setting $w(x)=\sqrt{1-x^2}$.
This weight function emphasizes the importance of the central part of the
interval $[-1,1]$ with respect to a neighborhood of the endpoints.
The corresponding recurrence coefficients are
\begin{equation}\label{cheby2}
\alpha_k=0, \ k=0,1,\ldots, \quad
\beta_k=\frac{1}{4}, \ k=1,2,\ldots, \quad
\beta_0=\frac{\pi}{2}.
\end{equation}

An image can be represented as the discretization of an intensity function
$f(x,y)$ on the continuous domain $[-1,1]^2:=[-1,1]\times[-1,1]$.
In such case, one may define bivariate polynomials as products of univariate
orthogonal polynomials
$$
P_{k,\ell}(x,y) = p_k(x) p_\ell(y).
$$
Such polynomials are themselves orthonormal, in the sense that
$$
\langle P_{k,\ell}, P_{r,s} \rangle =
\int_{-1}^1 \int_{-1}^1 P_{k,\ell}(x,y) P_{r,s}(x,y) w(x)w(y)\,dx\,dy 
= \delta_{kr} \delta_{\ell s},
$$
where $\delta_{ij}$ is the Kronecker symbol, that takes the value 1 when $i=j$
and 0 when $i\neq j$.

The orthogonal weighted moments of order $q$ of an intensity function $f(x,y)$ 
are defined as
\begin{equation}\label{Lmom_cont}
\mu_{i,q-i} = \int_{-1}^{1} \int_{-1}^{1} p_i(x) p_{q-i}(y) f(x,y)
w(x)w(y)\,dx\,dy,
\quad i=0,1,\ldots,q.
\end{equation}
The orthogonality of the basis implies that there is no redundancy or
overlapping of information between moments.
The knowledge of the moments, under suitable regularity conditions on the
intensity function, allows one to approximate $f(x,y)$ in the least squares
sense~\cite{riv81} in the form
\begin{equation}\label{LSapprox}
f_{mn}(x,y) = \sum_{i=0}^m \sum_{j=0}^n \mu_{ij}\, p_i(x) p_j(y).
\end{equation}
This means that
$$
\|f-f_{mn}\|_2 = \min_{p\in\Pi_{m,x}\otimes \Pi_{n,y}} \|f-p\|_2,
$$
where $\Pi_{m,x}$ is the vector space of polynomials in the variable $x$ having
degree lesser or equal than $m$, and $\otimes$ denotes the tensor product of
linear spaces.

Orthogonal moments \eqref{Lmom_cont} must be computed by a suitably accurate
quadrature formula. Given the large number of pixels usually available in an
image, we employed a simple Cartesian product rule based on Simpson's composite
formula \cite{sb91}.
Assuming the image has size $M\times N$, with both $M$ and $N$ odd, we denote
by
\begin{equation}\label{coors}
x_{k}=\frac{2k-(M+1)}{M-1}, \qquad y_{\ell}=\frac{2\ell-(N+1)}{N-1},
\end{equation}
($k=1,\ldots,M$ and $\ell=1,\ldots,N$)
the normalized pixel coordinates in the interval $[-1,1]$ and by
$f(x_k,y_\ell)$ the intensity of the pixel at row $k$ and column $\ell$.
Then, the moments are approximated as follows
\begin{equation}\label{simpson}
\mu_{ij} \approx K \sum_{k=1}^M \sum_{\ell=1}^N 
c_{M,k} c_{N,\ell}\, p_i(x_k) p_j(y_\ell) f(x_k,y_\ell),
\end{equation}
where $K=4/(9MN)$ and the quadrature weights are defined by
$$
c_{M,i} = \begin{cases}
1, \quad & i=1,M, \\
4, \quad & i=2,4,\ldots,M-1, \\
2, \quad & i=3,5,\ldots,M-2.
\end{cases}
$$

In some situations, it may be preferable to substitute \eqref{innprod} by the
following discrete inner product 
\begin{equation}\label{dinnprod}
\langle f, g \rangle = \sum_{k=1}^M f(t_k) g(t_k),
\end{equation}
where $\{t_k\}$ is a discretization of a given interval $[a,b]$.
If the interval $[0,M-1]$ is discretized by the points $t_k=k-1$,
$k=1,\ldots,M$, the polynomials $p^{(M)}_i(x)$, orthonormal on
$[0,M-1]$ with respect to \eqref{dinnprod}, are generally referred to as
discrete Chebyshev polynomials.
They can be obtained by substituting in \eqref{recurrence} the recurrence
coefficients 
$$
\alpha_k=\frac{M-1}{2}, \ k=0,1,\ldots, \quad
\beta_k=\frac{M^2}{4}\,\frac{1-\left(\frac{k}{M}\right)^2}{4-\frac{1}{k^2}},
\ k=1,2,\ldots, \quad
\beta_0=M.
$$
In such a case, all the theory is unchanged, except the continuous moments 
\eqref{Lmom_cont} are substituted by the discrete moments
\begin{equation}\label{discrmoms}
\mu_{i,q-i} = \sum_{k=1}^M \sum_{\ell=1}^N f_{k\ell} \, p^{(M)}_i(k-1) \,
p^{(N)}_{q-i}(\ell-1), \quad i=0,1,\ldots,q,
\end{equation}
where $f_{k\ell}$ is the intensity of the pixel at position $(k,\ell)$.

\section{Computation of orthogonal moments}\label{sec:comp}

The generation of orthogonal moments for an image can be a time consuming
process.
Indeed, formula \eqref{simpson} must be computed for each required moment.
The summation involves order $MN$ floating-point operations, once both
the polynomial bases have been evaluated at the points \eqref{coors}.
If the image is square ($M=N$) and the user sets $m=n$, only one basis
suffices, nevertheless the moments computation speed is strongly influenced by
the efficiency in the polynomial bases evaluation.
Using recurrence relations \eqref{recurrence}, the number of operations for the
evaluation of each basis is linear in the number of pixels on each side of the
image.

Additionally, it is essential to guarantee a sufficient accuracy in the
computation, since errors in the approximation of the moments may affect their
discriminative ability.
In some previous work, the moments have been approximated by considering in
\eqref{simpson} constant quadrature weights
\begin{equation}\label{badquad}
c_{M,k}=c_{N,\ell}=1 \quad \text{and} \quad K=\frac{(2i-1)(2j-1)}{MN}.
\end{equation}
This simplification greatly degrades the quality of the approximation of the
moments, without any advantage in terms of complexity.

There are various possible definitions of orthogonal polynomials.
Rodrigues' formula expresses them as derivatives of simple polynomials, which
is not useful in this setting.
Some authors adopted explicit expansions, like the following which holds for
Legendre polynomials,
\begin{equation}\label{badexp}
p_i(x) = 2^i \sum_{k=0}^i \binom{i}{k} \binom{\frac{i+k-1}{2}}{i} x^k.
\end{equation}
Such representations should be avoided because, as it is well known, the
canonical polynomial basis $\{1,x,x^2,\ldots,x^n\}$ is extremely
ill-conditioned~\cite{riv81}, in the sense that small perturbations in the
coefficients produce large errors in the polynomial evaluation.
Moreover, the coefficients in \eqref{badexp} have alternating signs, causing
cancellation and error propagation when the computation is performed on a
computer, and their implementation requires a large complexity, compared to
other approaches.

The Stieltjes procedure \eqref{recurrence}, despite being an implicit
definition, is the most effective from the computational point of view, as we
will show in the following.

\begin{algorithm}
\caption{[$\alphab$,$\betab$]=\opcoef(type,n,M) recurrence coefficients of
orthogonal polynomials}
\label{algo:opcoef}
\begin{algorithmic}[1]
\REQUIRE orthogonal polynomial type, maximum degree $n$, $M$ number of
discretization points (only for discrete Chebyshev polynomials)
\ENSURE $\alphab,\betab\in\R^{n+1}$ vectors of coefficients for the
orthogonal polynomials
\STATE $\bm{v}=(1,2,\ldots,n+1)$
\IF{type='Legendre'}
	\STATE $\alphab=(0,\ldots,0)$
	\STATE $\betab=\big(2,\bm{v}^2/(4\bm{v}^2-1)\big)$ (componentwise)
		\label{opcoef:line4}
\ELSIF{type='Chebyshev second kind'}
	\STATE $\alphab=(0,\ldots,0)$
	\STATE $\betab=\big(\frac{\pi}{2},\frac{1}{4},\ldots,\frac{1}{4}\big)$
\ELSIF{type='discrete Chebyshev'}
	\STATE $\alphab=\frac{M-1}{2}(1,\ldots,1)$
	\STATE $\betab=\big(M,\frac{M^2}{4}(1-(\bm{v}/M)^2)/(4-1/\bm{v}^2)\big)$
		(componentwise) \label{opcoef:line10}
\ENDIF
\end{algorithmic}
\end{algorithm}

Algorithm~\ref{algo:opcoef} implements a simple function \opcoef\ which
returns the recurrence coefficients of a few families of orthogonal
polynomials. At the moment, only the families mentioned in
Section~\ref{sec:legmom} are supported, but the function can be easily
extended.
The operations involving vectors (displayed as bold lower
case letters) in lines~\ref{opcoef:line4} and~\ref{opcoef:line10} of
Algorithm~\ref{algo:opcoef} are to be performed component by component.
This implies, in Matlab, the use of the array operators \texttt{.*},
\texttt{./}, and \texttt{.\^}, in place of the usual matrix operators.

\begin{algorithm}
\caption{$P$=\opevmat($\alphab$,$\betab$,$\bm{x}$) evaluation of orthogonal
polynomials}
\label{algo:opevmat}
\begin{algorithmic}[1]
\REQUIRE $\alphab,\betab\in\R^{n+1}$ vectors of coefficients for the
orthogonal polynomials, $\bm{x}\in\R^m$ vector of evaluation points
\ENSURE $P\in\R^{m\times(n+1)}$ matrix whose column $P_{:,j}$ is the evaluation
of the $(j-1)$th orthogonal polynomial at the components of $\bm{x}$
\STATE $P_{:,1}=(1,\ldots,1)^T/\sqrt{\beta_1}$
\STATE $P_{:,2}=(\bm{x}-\alpha_1) * P_{:,1}$ (componentwise)
	\label{opevmat:line2}
\STATE $P_{:,2}=P_{:,2}/\sqrt{\beta_2}$
\FOR{$k=2,\ldots,n$}
	\STATE $P_{:,k+1}=(\bm{x}-\alpha_k) *
		P_{:,k}-\sqrt{\beta_k}\,P_{:,k-1}$ (componentwise)
		\label{opevmat:line5}
	\STATE $P_{:,k+1}=P_{:,k+1}/\sqrt{\beta_{k+1}}$
\ENDFOR
\end{algorithmic}
\end{algorithm}

Once the recursion coefficients have been assigned, the evaluation of an
orthogonal polynomials basis at a vector of points is performed by the function
\opevmat\ in Algorithm~\ref{algo:opevmat}, which employs the Stieltjes
procedure.
Again, vector multiplications in lines~\ref{opevmat:line2}
and~\ref{opevmat:line5} are intended component by component.

\begin{algorithm}
\caption{$M$=\legmoms($F$,$q$) computation of Legendre moments}
\label{algo:legmoms}
\begin{algorithmic}[1]
\REQUIRE $F\in\R^{M\times N}$ grayscale digital image, $q$ order of the moments
to be computed
\ENSURE $M\in\R^{(q+1)\times(q+1)}$ matrix of Legendre moments
\STATE \textbf{if} $M$ is even \textbf{then} $M=M-1$ \textbf{end if}
\STATE \textbf{if} $N$ is even \textbf{then} $N=N-1$ \textbf{end if}
\STATE $h_x=2/M$, $\bm{x}=(-1,-1+h_x,-1+2h_x,\ldots,1)^T$
\STATE $h_y=2/N$, $\bm{y}=(-1,-1+h_y,-1+2h_y,\ldots,1)^T$
\STATE [$\alphab$,$\betab$]=\opcoef('Legendre',q)
\STATE $P_1$=\opevmat($\alphab$,$\betab$,$\bm{x}$)
\STATE $P_2$=\opevmat($\alphab$,$\betab$,$\bm{y}$)
\FOR{$i=1,\ldots,q+1$}
	\FOR{$j=1,\ldots,q+1$}
		\STATE compute $M_{i-1,j-1}$ by applying Simpson's rule
		\eqref{simpson} to integral \eqref{Lmom_cont}
	\ENDFOR
\ENDFOR
\end{algorithmic}
\end{algorithm}

Algorithm~\ref{algo:legmoms} describes the function \legmoms\ which
approximates the orthogonal moments \eqref{Lmom_cont} with respect to Legendre
polynomials by the Cartesian product Simpson's rule \eqref{simpson}.
The approximation of the orthogonal moments with respect to Chebyshev
polynomials of the second kind is performed by a very similar function
\chebtwomoms, which we do not report here for the sake of brevity.

\begin{algorithm}
\caption{$M$=\dchebmoms($F$,$q$) comp. of discrete Chebyshev moments}
\label{algo:dchebmoms}
\begin{algorithmic}[1]
\REQUIRE $F\in\R^{M\times N}$ digital grayscale image, $q$ order of the moments
to be computed
\ENSURE $M\in\R^{(q+1)\times(q+1)}$ matrix of Legendre moments
\STATE $\bm{x}=(0,1,\ldots,M-1)^T$
\STATE $\bm{y}=(0,1,\ldots,N-1)^T$
\STATE [$\alphab_1$,$\betab_1$]=\opcoef('discrete Chebyshev',q,M)
\STATE [$\alphab_2$,$\betab_2$]=\opcoef('discrete Chebyshev',q,N)
\STATE $P_1$=\opevmat($\alphab_1$,$\betab_1$,$\bm{x}$)
\STATE $P_2$=\opevmat($\alphab_2$,$\betab_2$,$\bm{y}$)
\STATE $M=P_1^T*F*P_2$ (matrix product)
\end{algorithmic}
\end{algorithm}

Finally, the exact computation of the discrete orthogonal moments
\eqref{discrmoms} is described in the function \dchebmoms\ of
Algorithm~\ref{algo:dchebmoms}.

All the above algorithms have been implemented in the Matlab programming
language and are available from the authors' web sites; see the
\texttt{orthomoms} package at \url{http://bugs.unica.it/cana/software}.
The software have been verified to be compatible with Octave.

\section{Experimental results}\label{sec:results}

In this section we present two sets of numerical experiments carried out in
order to verify the effectiveness of the proposed algorithm.
The experiments were performed on a dual Xeon CPU E5-2620 system (12 cores),
running the Debian GNU/Linux operating system and Matlab 9.2.

In the first experiment, we investigate the reconstruction performance of the
above described approach for computing Legendre moments.
Such an experiment is useful to ascertain the correctness of the computation
and the advantages, both in terms of execution time and of accuracy, of the
recurrence relations with respect to the closed form representation of Legendre
polynomials.

In the second experiment, the descriptive ability of the computed moments is
tested in texture classification problems.

\subsection{Image reconstruction}

\begin{figure}[htb]
\begin{center}
\includegraphics[width=.55\textwidth]{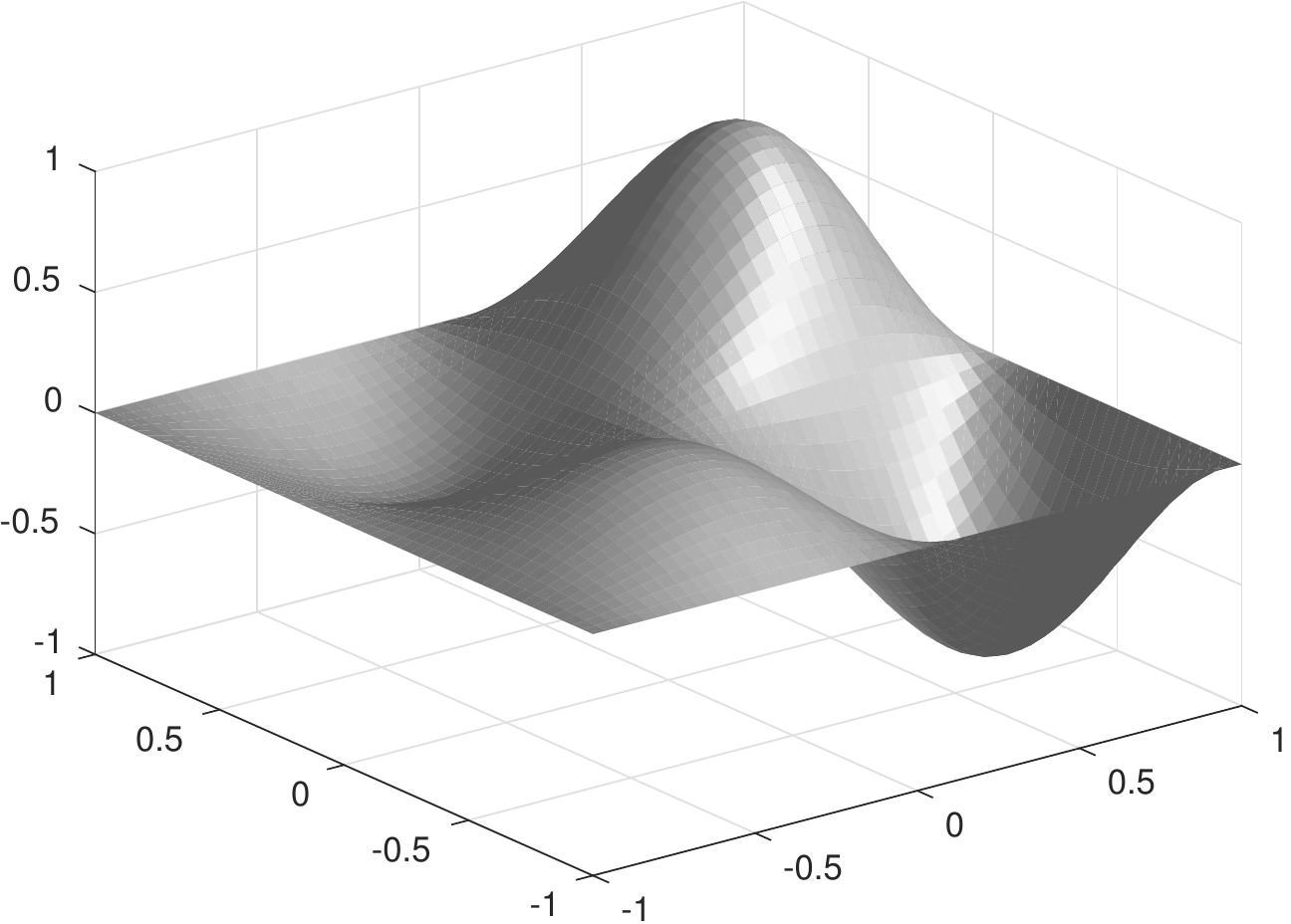}
\includegraphics[width=.40\textwidth]{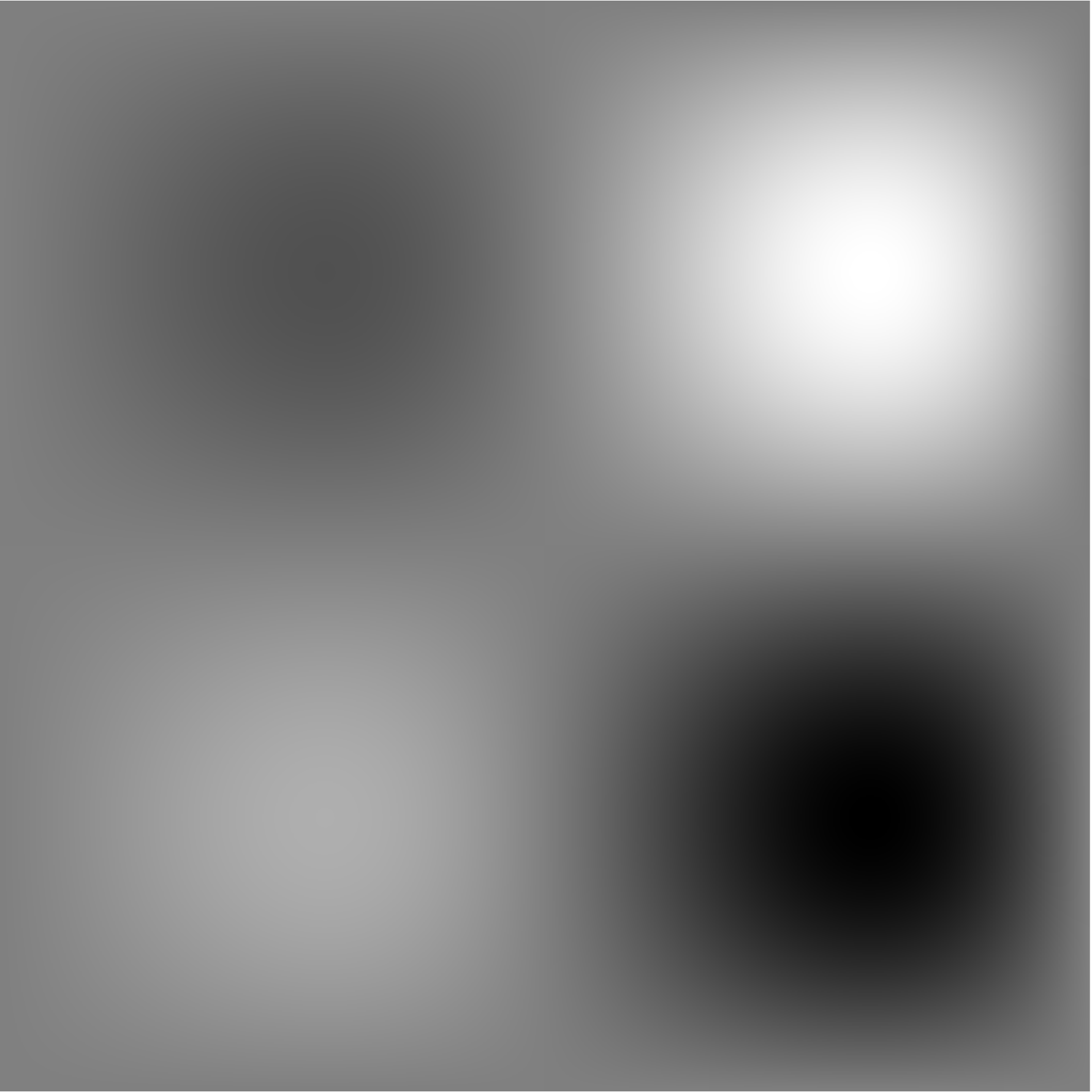}
\caption{Model intensity function and corresponding grayscale image.}
\label{fig:model}
\end{center}
\end{figure}

In order to verify the correctness of the computation on a synthetic dataset,
we consider the following smooth model intensity function
$$
f(x,y) = \frac{1}{2} \e^x \sin(\pi x) \sin(\pi y), \qquad (x,y)\in[-1,1]^2.
$$
The function is sampled on a regular grid of $1023\times 1023$ points on the
square $[-1,1]^2$, then its values are translated and rescaled so that they are
contained in the interval $[0,1]$.
The graph of the function and the resulting grayscale image $F$ are displayed
in Figure~\ref{fig:model}.

We first evaluate the moments $\mu_{ij}$, $i,j=0,\ldots,n$, by Simpson's
rule~\eqref{simpson}, where the orthogonal polynomials are computed using the
three-term recursion \eqref{recurrence}.
Then the least-squares approximation $f_{nn}(x,y)$ of the model intensity
function $f(x,y)$ is evaluated by formula~\eqref{LSapprox} on the
grid~\eqref{coors}, as well as the relative error 
$$
E_n(f)=\frac{\max_{k,\ell} |f(x_k,y_\ell)-f_{nn}(x_k,y_\ell)|}{\max_{k,\ell}
|f(x_k,y_\ell)|}.
$$
We perform different tests, letting the number of moments $n=5,10,\ldots,50$.

\begin{figure}[htb]
\begin{center}
\includegraphics[width=.49\textwidth]{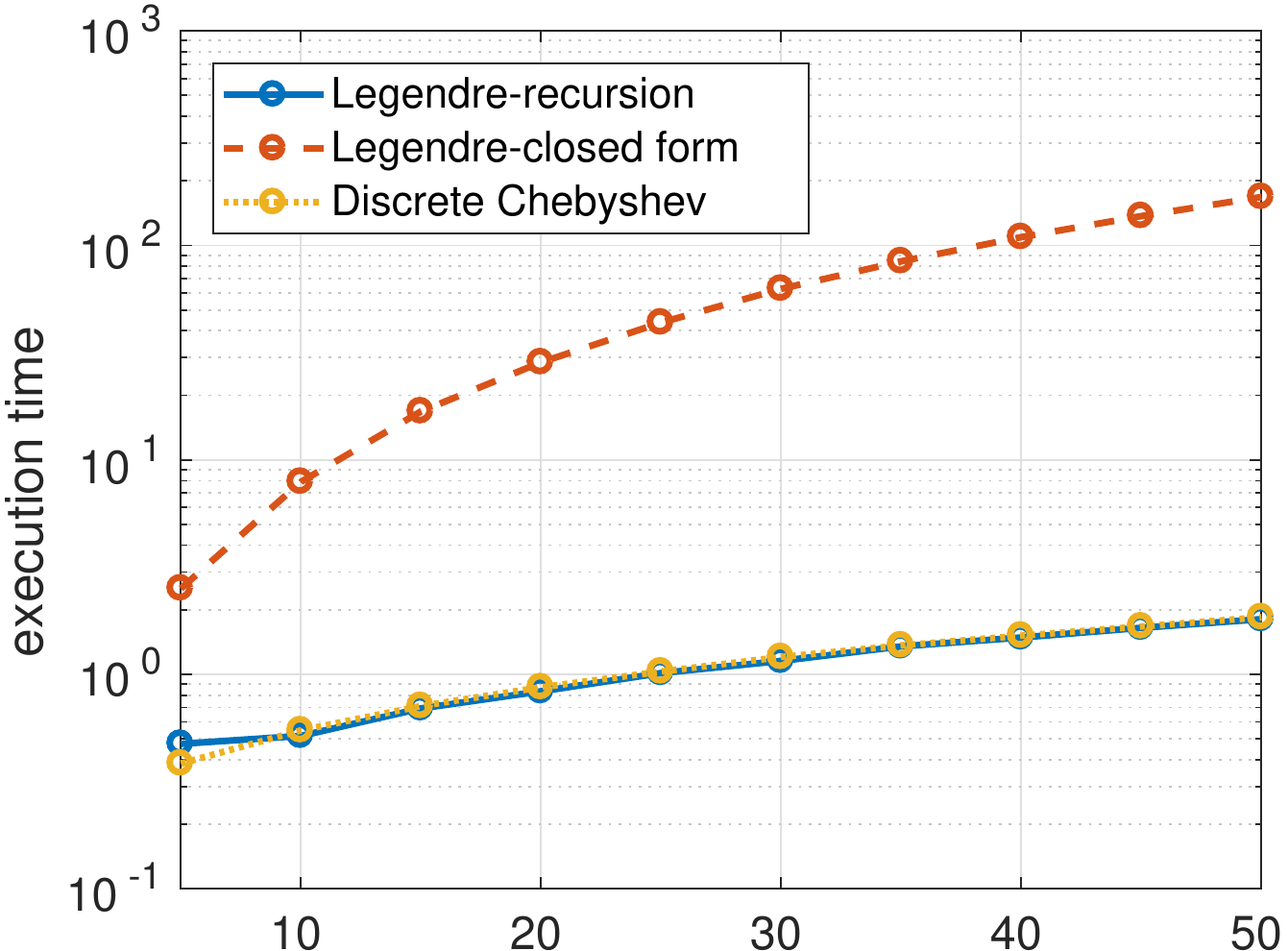}
\includegraphics[width=.49\textwidth]{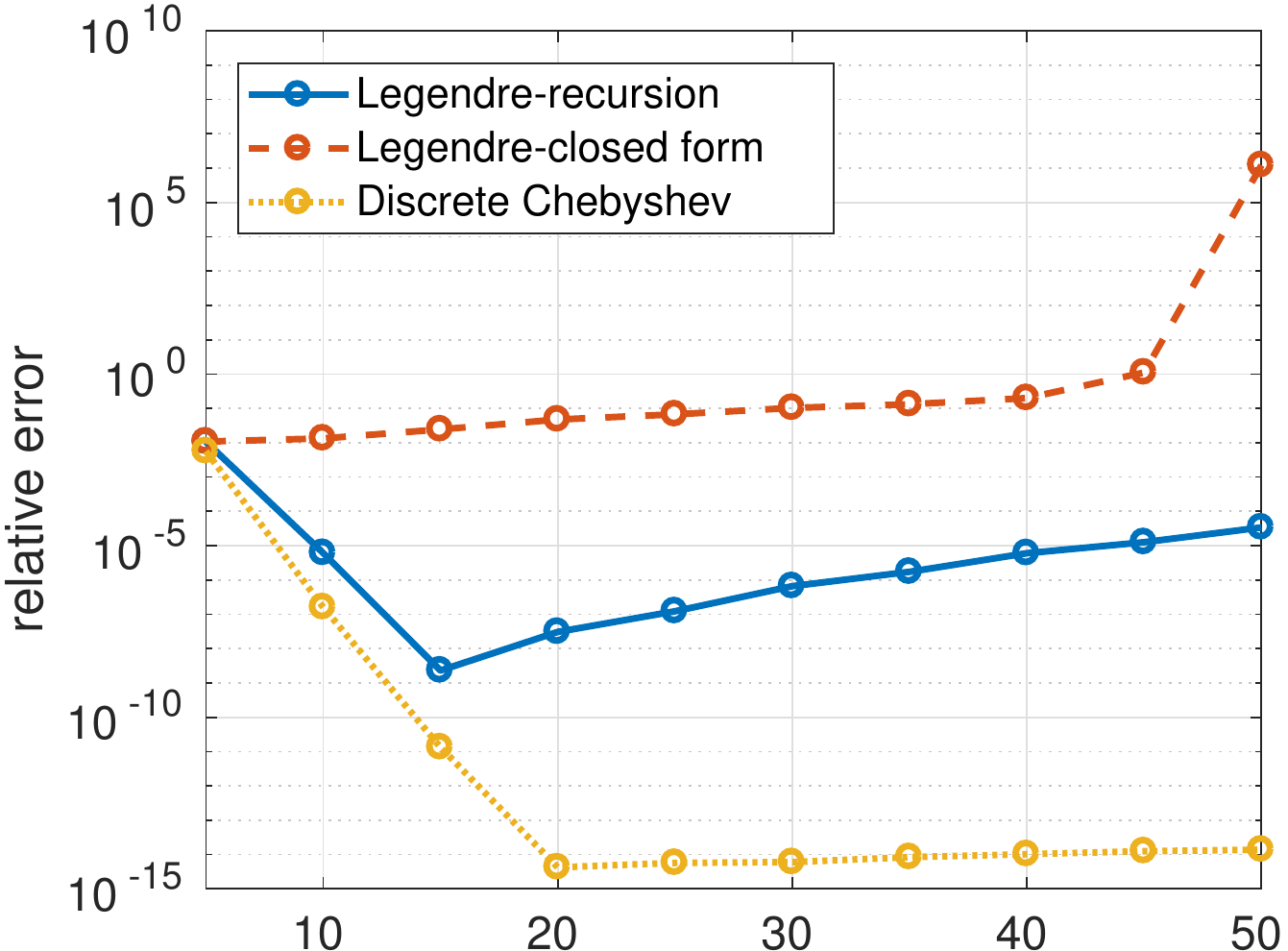}
\caption{Execution time and reconstruction error when the number of moments
takes the values $n=5,10,\ldots,50$.}
\label{fig:reconst}
\end{center}
\end{figure}

\begin{figure}[!htb]
\begin{center}
\includegraphics[width=.47\textwidth]{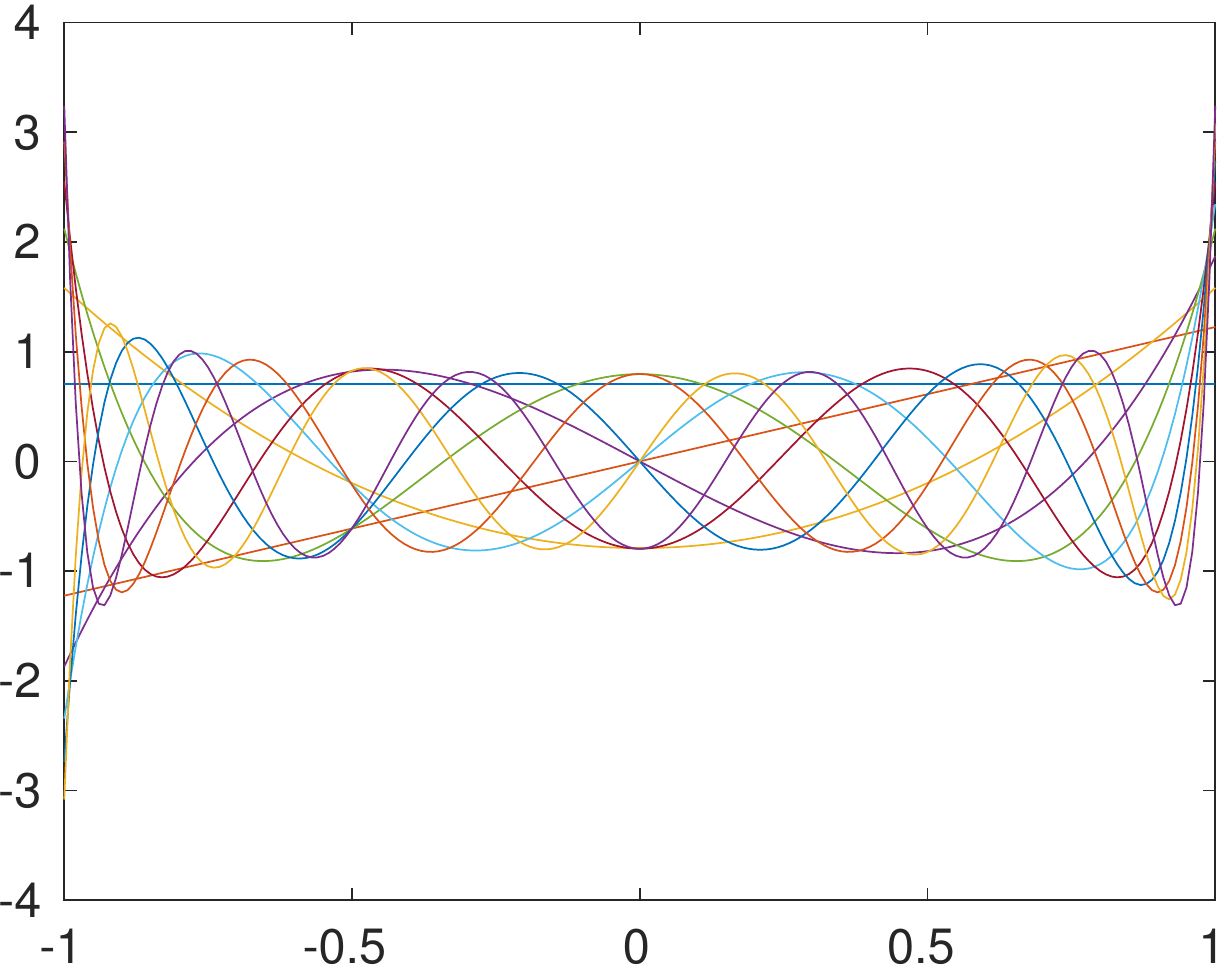}
\hfill
\includegraphics[width=.50\textwidth]{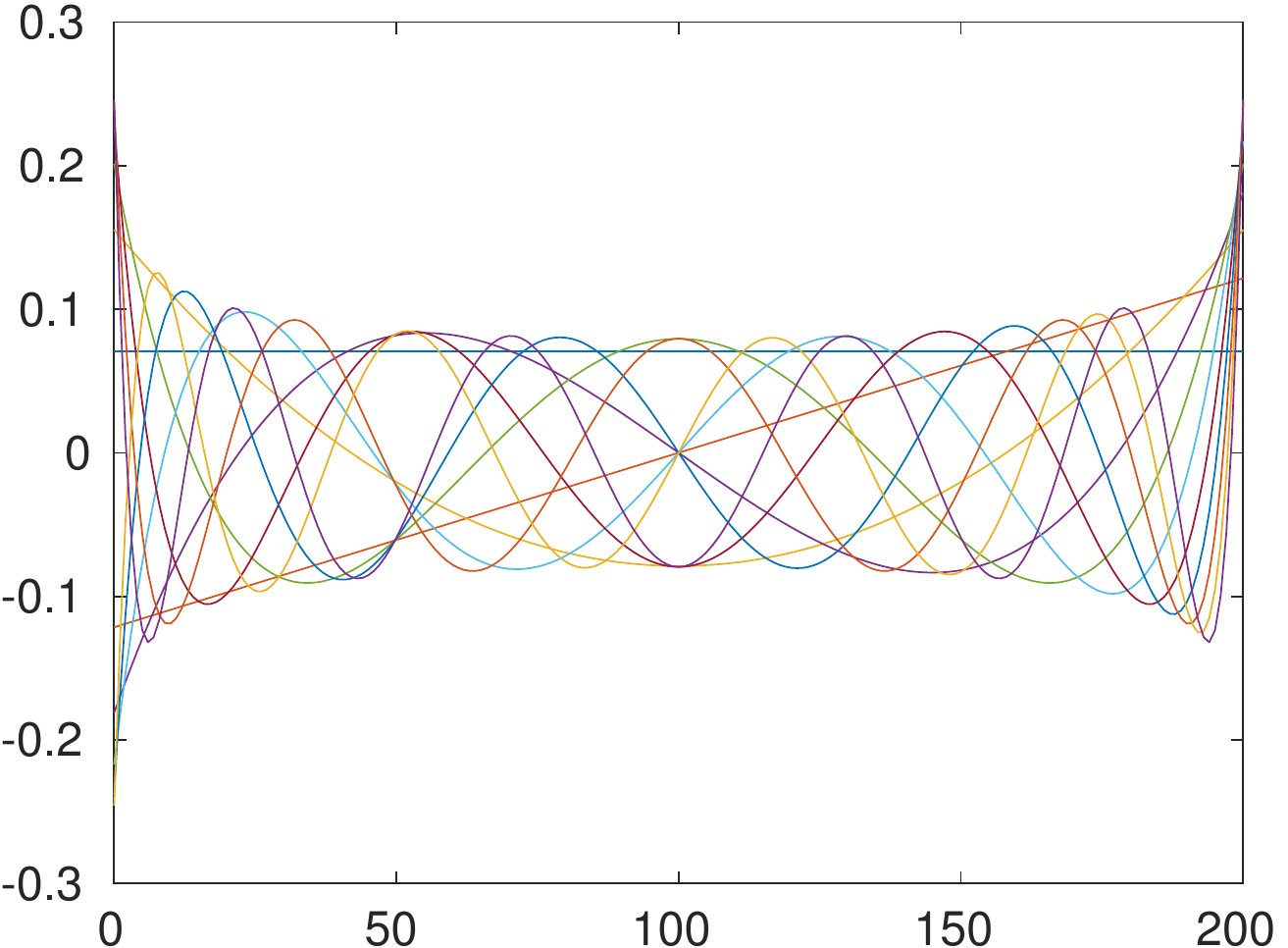}
\caption{Orthonormal Legendre polynomials (on the left) and discrete Chebyshev
polynomials (on the right) of degree $n=0,1,\ldots,10$.}
\label{fig:polys}
\end{center}
\end{figure}

In Figure~\ref{fig:reconst} we compare the execution time and the relative
error in the reconstruction, obtained by the above approach, to the results
obtained by evaluating the Legendre polynomials using the closed form
representation \eqref{badexp} and approximating the moments by employing in
\eqref{simpson} the quadrature weights \eqref{badquad}.

The graph on the left of Figure~\ref{fig:reconst} shows the great reduction in
the computational complexity of the approach based on recurrence relations.
The reconstruction error is displayed on the right.
The error deriving from the closed form representation increases with the
number of moments, because of error propagation and lack of precision in the
approximation of the integrals.
The computation performed by recursion formulas appears to be much more stable.
The accuracy improves up to 15 moments, then it degrades, but the relative
error stays below $10^{-4}$. This degradation is probably due to the fact that
Simpson's rule is not able to approximate with sufficient accuracy the integral
of fast oscillating functions, like large degree orthogonal polynomials are.

The same figure displays also the results obtained by discrete Chebyshev
polynomials and the corresponding discrete moments \eqref{discrmoms}.
The execution time is equivalent to Legendre polynomials, indeed the two
approaches use the same computational scheme with different recursion
coefficients. 
The accuracy improves up to 20 moments, where it reaches machine
accuracy.
This is somehow expected, as in this case the moments are computed exactly and
not obtained by approximating an integral by a quadrature formula.
Nevertheless, the fact that the accuracy stays at the same level when $n>20$
clearly testifies that error propagation is negligible when recursion formulas
\eqref{recurrence} are used.

The two families of orthogonal polynomials used in the experiments are depicted
in Figure~\ref{fig:polys}.

\subsection{Texture analysis}

\begin{figure}[!tb]
\begin{center}
\centerline{\includegraphics[width=0.98\textwidth]{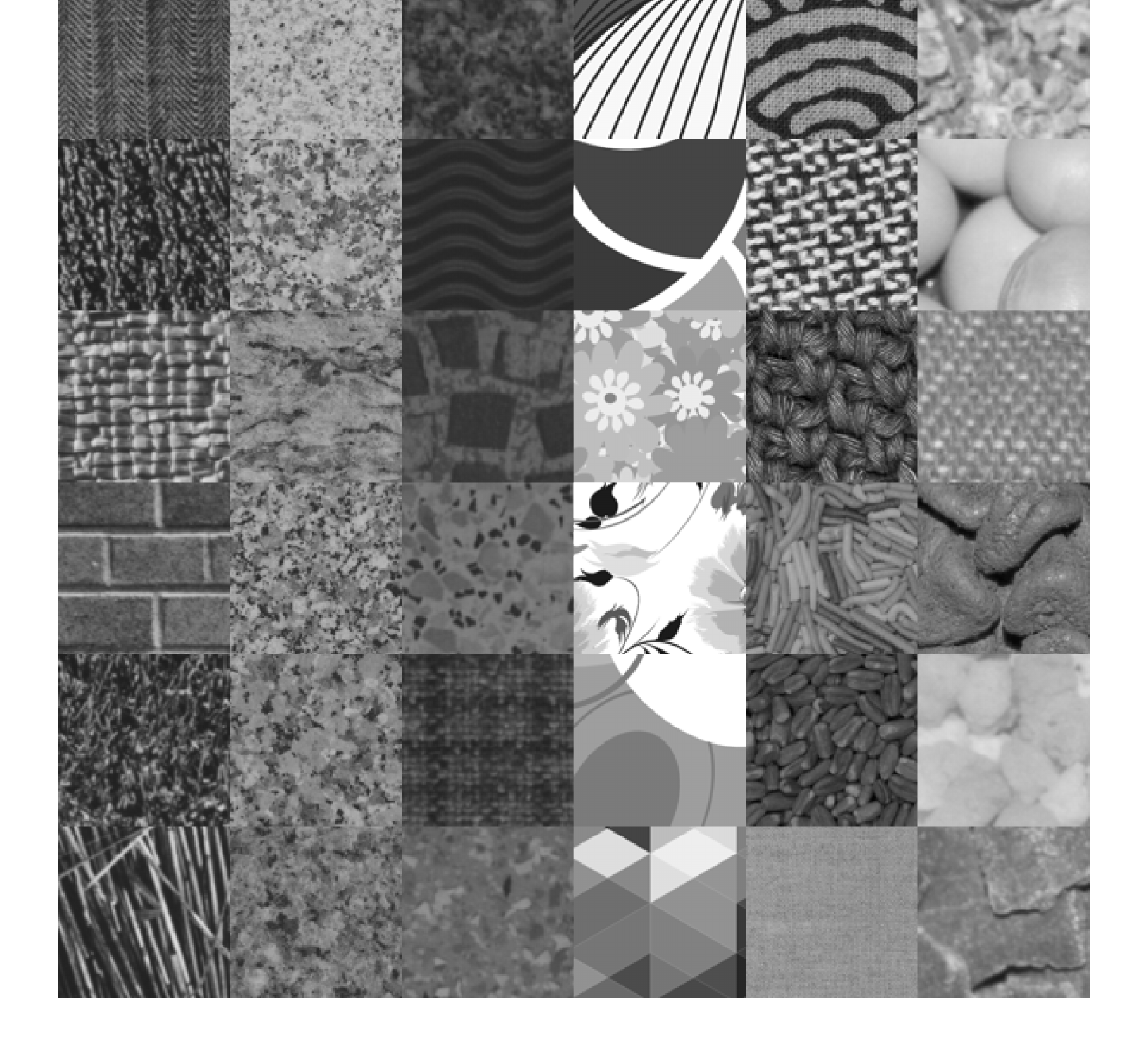}}
\caption{Each column displays images from the six texture databases used in the
experiments. Starting from the left: Brodatz, Mondial Marmi, Outex, Vectorial,
Kylberg Sintorn, and ALOT.}
\label{fig:textures}
\end{center}
\end{figure}

To evaluate and compare the performance in texture analysis of the analyzed
algorithms for computing orthogonal moments, we used six different databases of
texture images: Brodatz, Mondial Marmi, Outex, Vectorial, Kylberg Sintorn and
ALOT; see Figure~\ref{fig:textures}. 
They present various materials and textures representing different
classification problems, and include hardware-rotated images taken at 4--9
different orientations, making them the most suitable for our experiments.

The only database that contains software rotated images is Vectorial, a
collection of 20 artificial texture classes proposed by
Bianconi~\cite{Bianconi}.
As it can be guessed from the name, these images are not raster, thus the
rotation does not affect the image structure. They have been rotated with
angular steps of $10^\circ$ ($0^\circ,10^\circ,\ldots,90^\circ$) and then
converted into raster with a resolution of 300 dpi.
Finally, each image has been subdivided into 16 sub-images of size
$255 \times 255$, resulting in 16 samples per class and 3200 total images.
We did not consider other databases which include software rotated images,
since this operation may modify the original image structure, leading to
wrong results.

The other five databases present real texture images.
The first one belongs to one of the most
diffused collection of textures, that is, Brodatz's album. 
We did not use the whole collection, since the original images included
in the album are not rotated, but just a subset of 13 textures. This
subset has been proposed by Bianconi \cite{Bianconi}, who acquired
hardware-rotated images directly from the original book. The images have been
acquired with angular steps of $10^\circ$, like in the previous case.
Then, each image has been subdivided into 16 sub-images of size
$205 \times 205$, resulting in 2080 total samples. Brodatz's album is one of
the oldest texture database. Indeed, among the tested datasets it is the only
one composed by gray level images.

The remaining databases comprise color images, thus they have been converted
into gray scale before the feature extraction step.
Mondial Marmi is a free image database of
granite tiles for color and texture analysis, that includes 12 granite classes.
Each texture is captured in a 24-bit RGB image of size $544\times 544$ using
nine rotation angles ($0^\circ,5^\circ,10^\circ,15^\circ,30^\circ,
45^\circ,60^\circ,75^\circ,90^\circ$). The database presents 4
images for each class and for each angle, so the original images have been
subdivided into four non-overlapping sub-images of size $272 \times 272$ for a
total image count of 1780.

The Outex database is a collection of 320
textures (both macrotextures and microtextures) acquired with well defined
variations in terms of illumination, rotation, and spatial resolution. Each
texture is captured in a 24-bit RGB image of size $538\times 746$, using three
different simulated illuminations, six spatial resolutions ($100,120,300,360,
500$, and 600 dpi), and nine rotation angles ($0^\circ,5^\circ,10^\circ,
15^\circ,30^\circ,45^\circ,60^\circ,75^\circ,90^\circ$).
Hence the current texture database includes 51840 images. Given the
considerable size of this database, we focused our experiments on a smaller
subset, taking inspiration by the test suite proposed by the author of
the Outex database himself, called OUTEX00045, that uses 45 texture classes.
Following the instructions proposed in Outex, the original images have been
divided in 20 non overlapping sub-images with size $128 \times 128$, for a
total count of 8100 images.

The Kylberg Sintorn database is a collection of 25 textural
classes of materials, such as fabric, grains, sugar, rice, etc. As in the other
databases, the images are provided with nine rotation angles, but in this case
the images have been rotated with angular steps of $40^\circ$
($0^\circ,40^\circ,\ldots,320^\circ$). The original images (one for each class)
are 24-bit RGB with a resolution of $5184 \times 3456$ pixels, but they have
been provided also in small subsets for texture classification, presenting 400
images for each angle and thus 16 samples per class. The final dataset
contains, therefore, 3600 images.

The last database, ALOT, is very different from the previous ones. Indeed, it
presents images acquired with just four rotation angles 
with steps of $60^\circ$ ($0^\circ,60^\circ,120^\circ,180^\circ$).
The original image database contains 250 textures, each one with
100 images obtained under different illumination conditions. For our
experiments we considered a subset of the original dataset containing 
80 textures. The original images have been divided into 16 sub-images of
size $181\times 181$, for a total of 5120 images.

To compare the recurrent formulation of orthogonal moments to the one based on
a closed form expansion, we performed a set of image classification
experiments, evaluating both performance and robustness of the proposed
descriptors against image rotation.
For each dataset, we
executed 100 experiments; for each of them, both the training and the test sets
are represented by half of the original samples. Each dataset was
divided by a stratified sampling, which guarantees that each class is
properly represented in both the training and the testing set. The
classification performances have been evaluated by the accuracy index
\begin{equation}\label{accuracy}
\mathcal{A} = \frac{TP + TN}{TP + TN + FP + FN},
\end{equation}
where TP (TN) and FP (FN) represent the number of true positive (negative)
classifications, and false positive (negative) classifications, respectively.
This index gives a meaningful indication of the performance, since it considers
each class of equal importance.

The classification accuracy has been estimated by a k-Nearest Neighbor
(k-NN) classifier, with $k = 1$, computed using the Euclidean distance.
The k-NN has been preferred to a more complex classifier in order to document
the effectiveness of the moments as features descriptors, rather than assessing
the performance of the classifier itself. 
To better study the effects of image rotation, the classifier is always trained
with features extracted from images acquired at orientation $0^\circ$ and then
tested with feature extracted from images acquired at other orientations.

\begin{figure}[!ht]
\begin{center}
\centerline{\includegraphics[width=.48\textwidth]{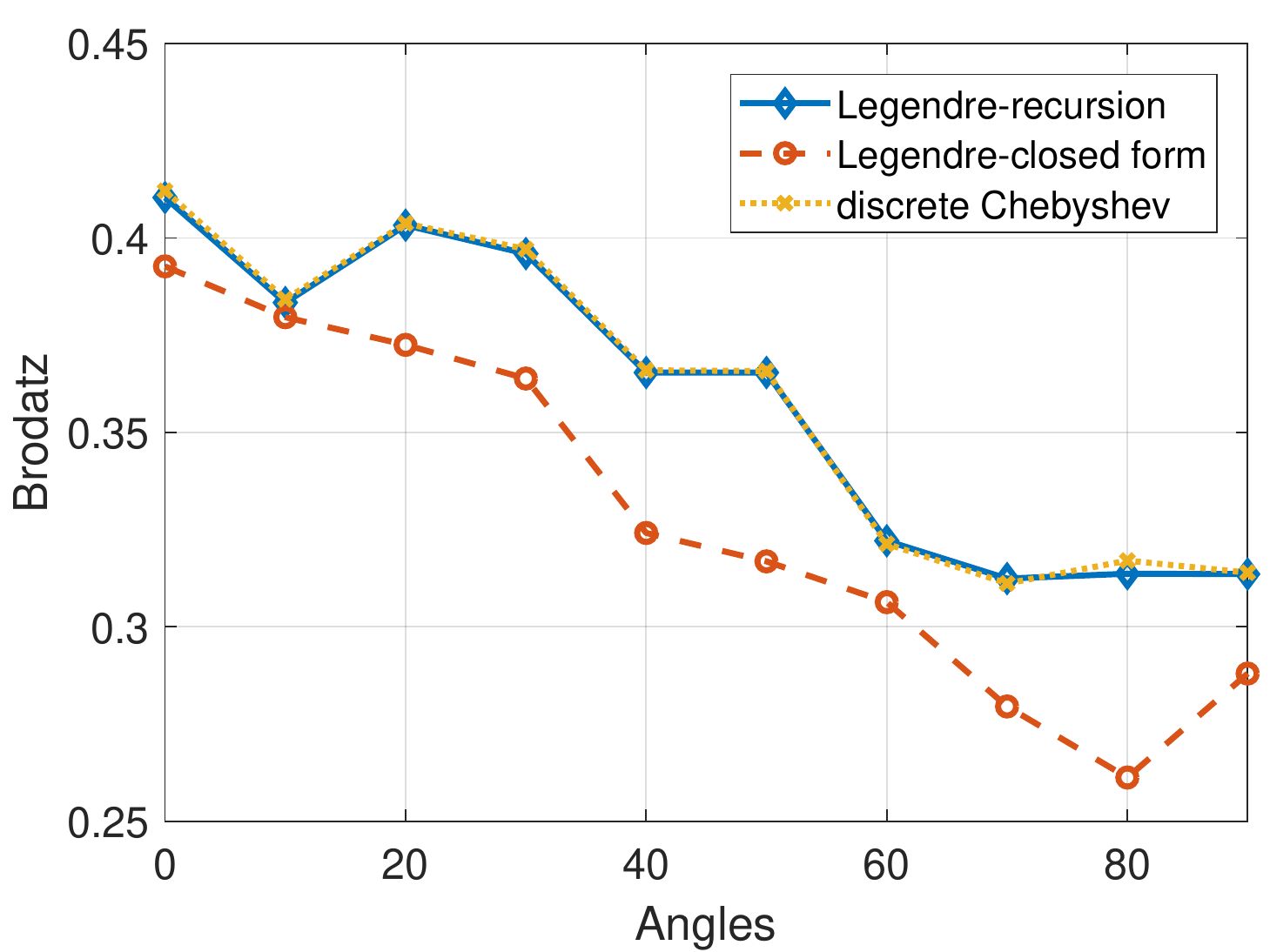} \hfill \includegraphics[width=.48\textwidth]{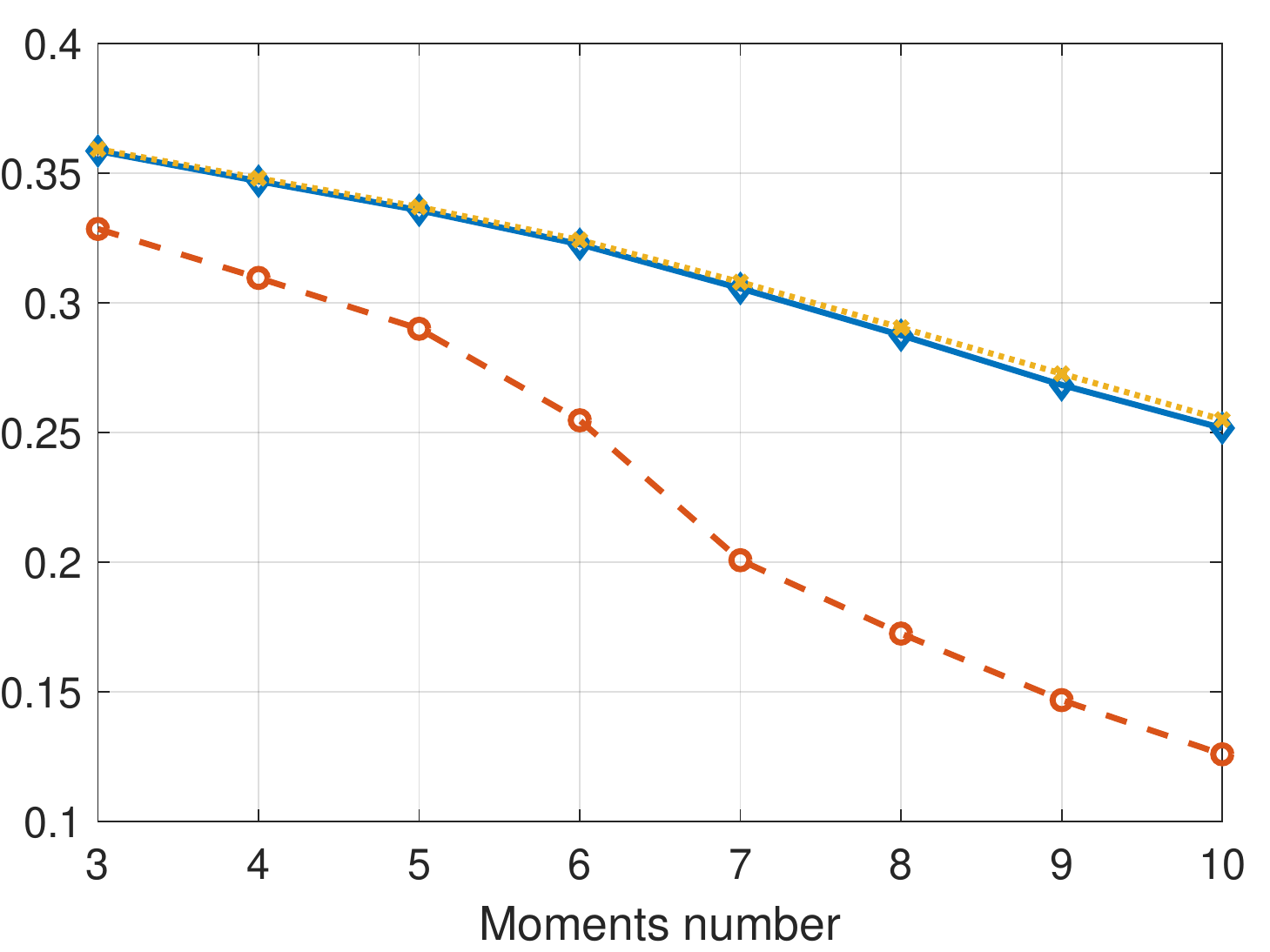}}
\bigskip
\centerline{\includegraphics[width=.48\textwidth]{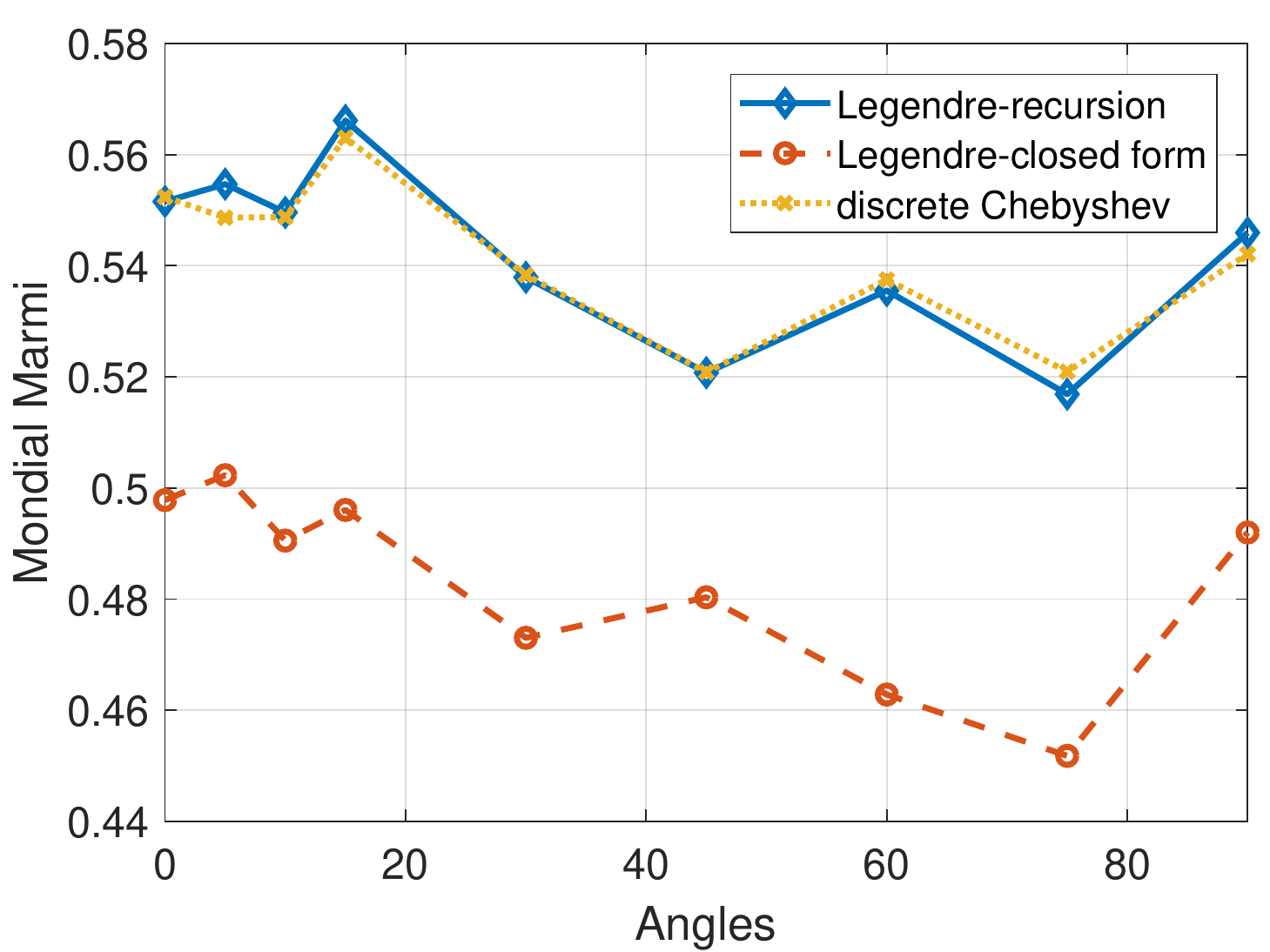} \hfill \includegraphics[width=.48\textwidth]{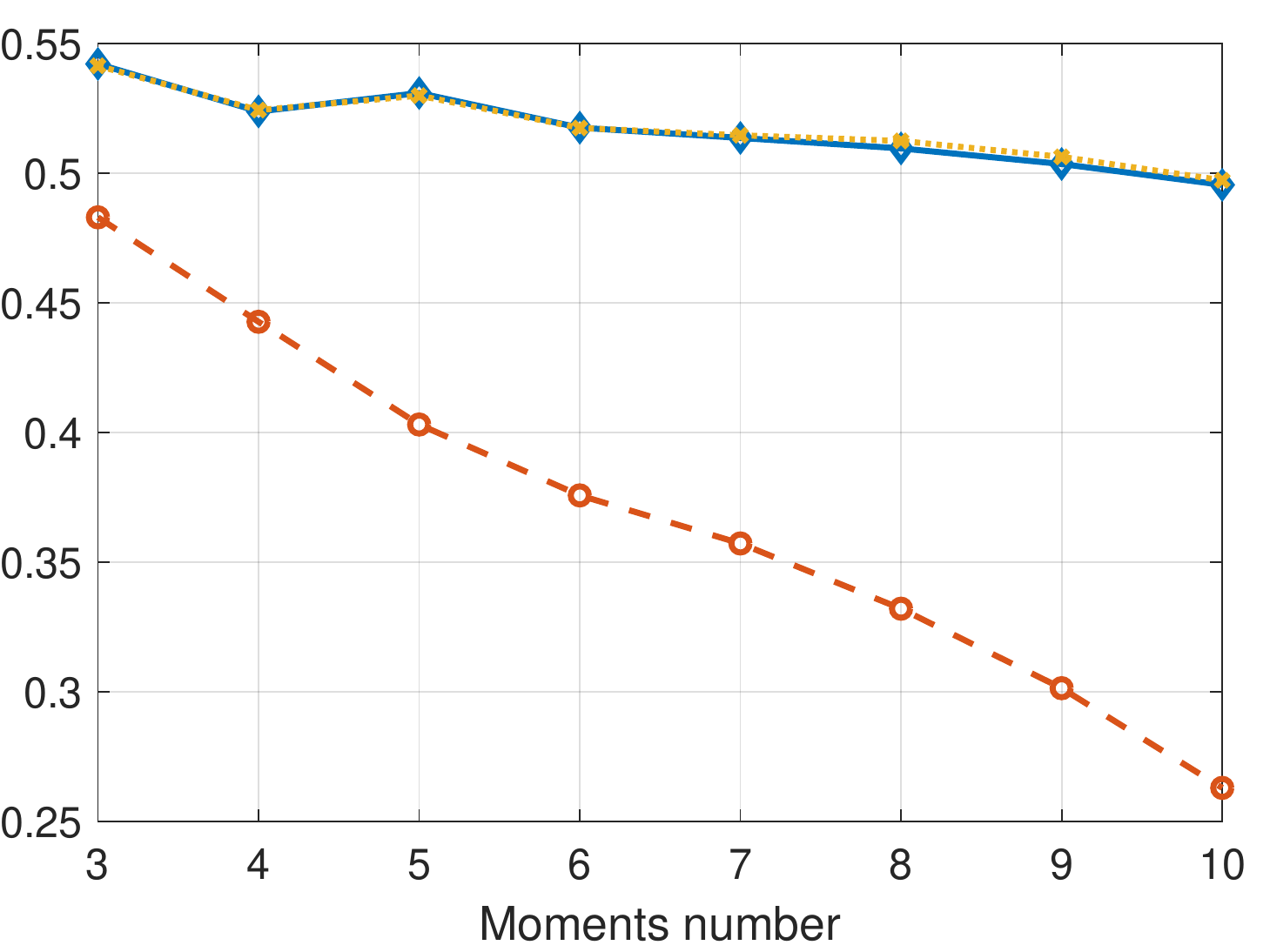}}
\bigskip
\centerline{\includegraphics[width=.48\textwidth]{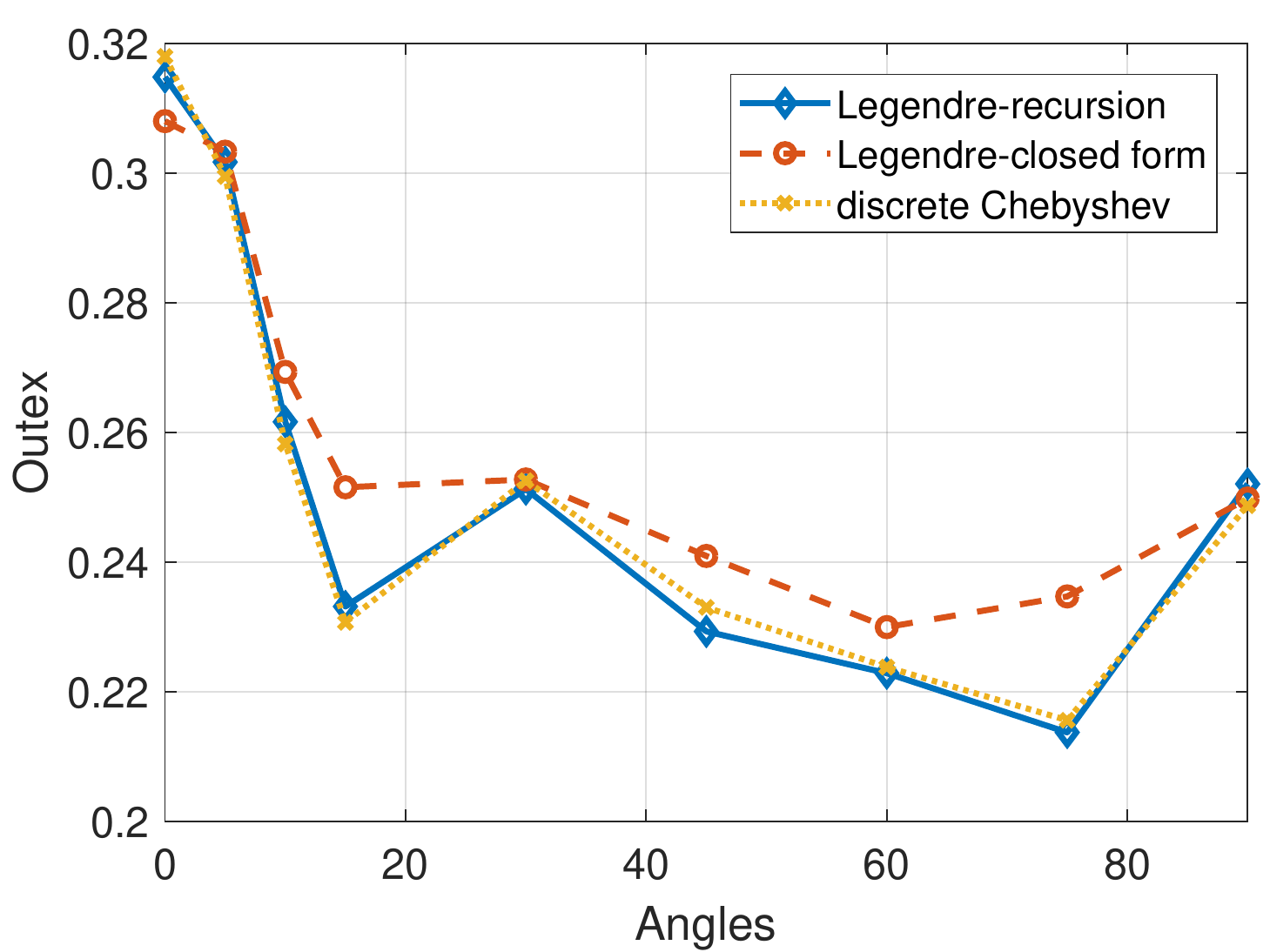} \hfill \includegraphics[width=.48\textwidth]{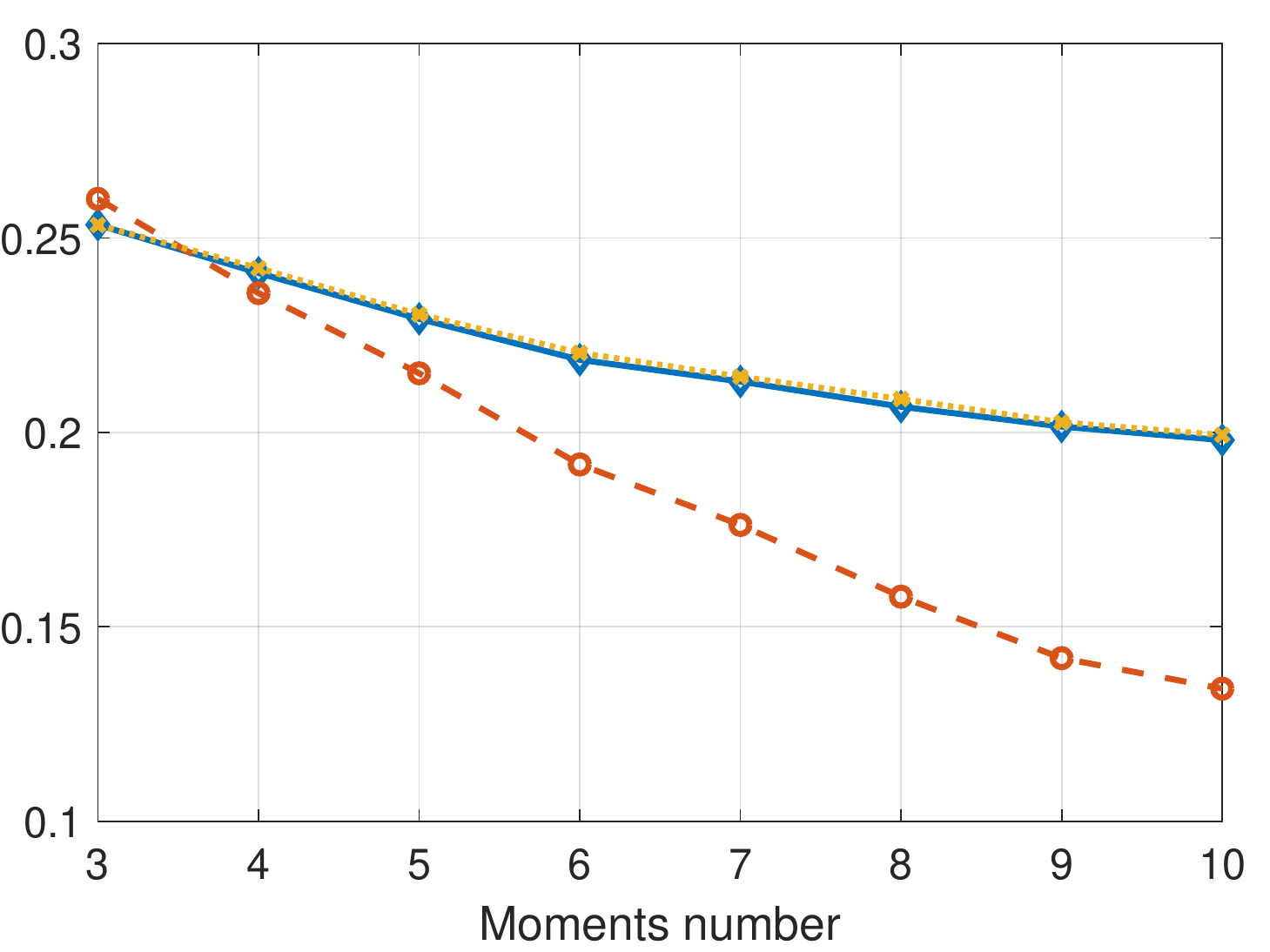}}
\caption{Brodatz, Mondial Marmi, and Outex databases: accuracy index
\eqref{accuracy} versus image orientation (left) and number of moments (right).}
\label{fig:performances}
\end{center}
\end{figure}

\begin{figure}[!ht]
\begin{center}
\centerline{\includegraphics[width=.48\textwidth]{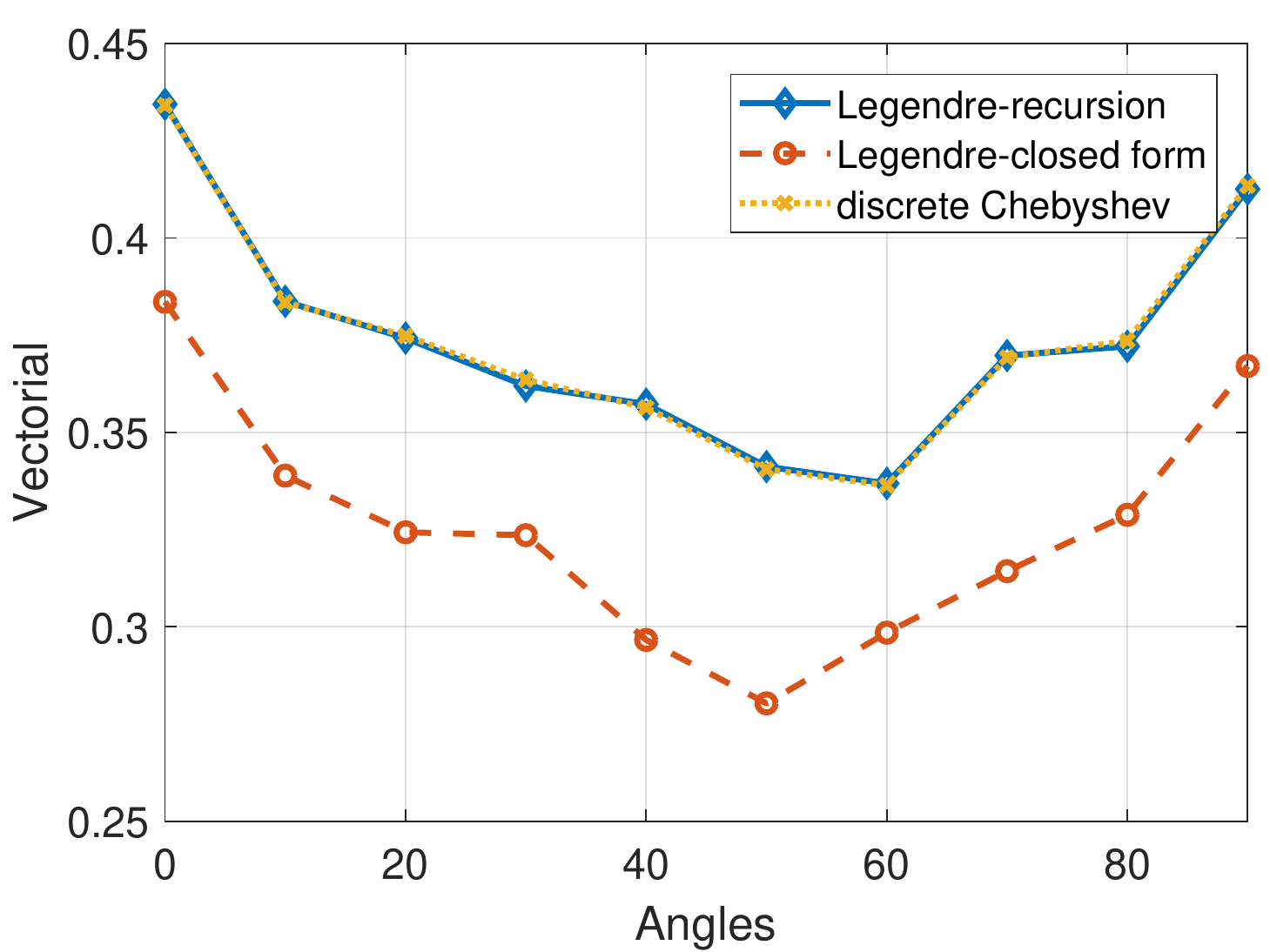} \hfill \includegraphics[width=.48\textwidth]{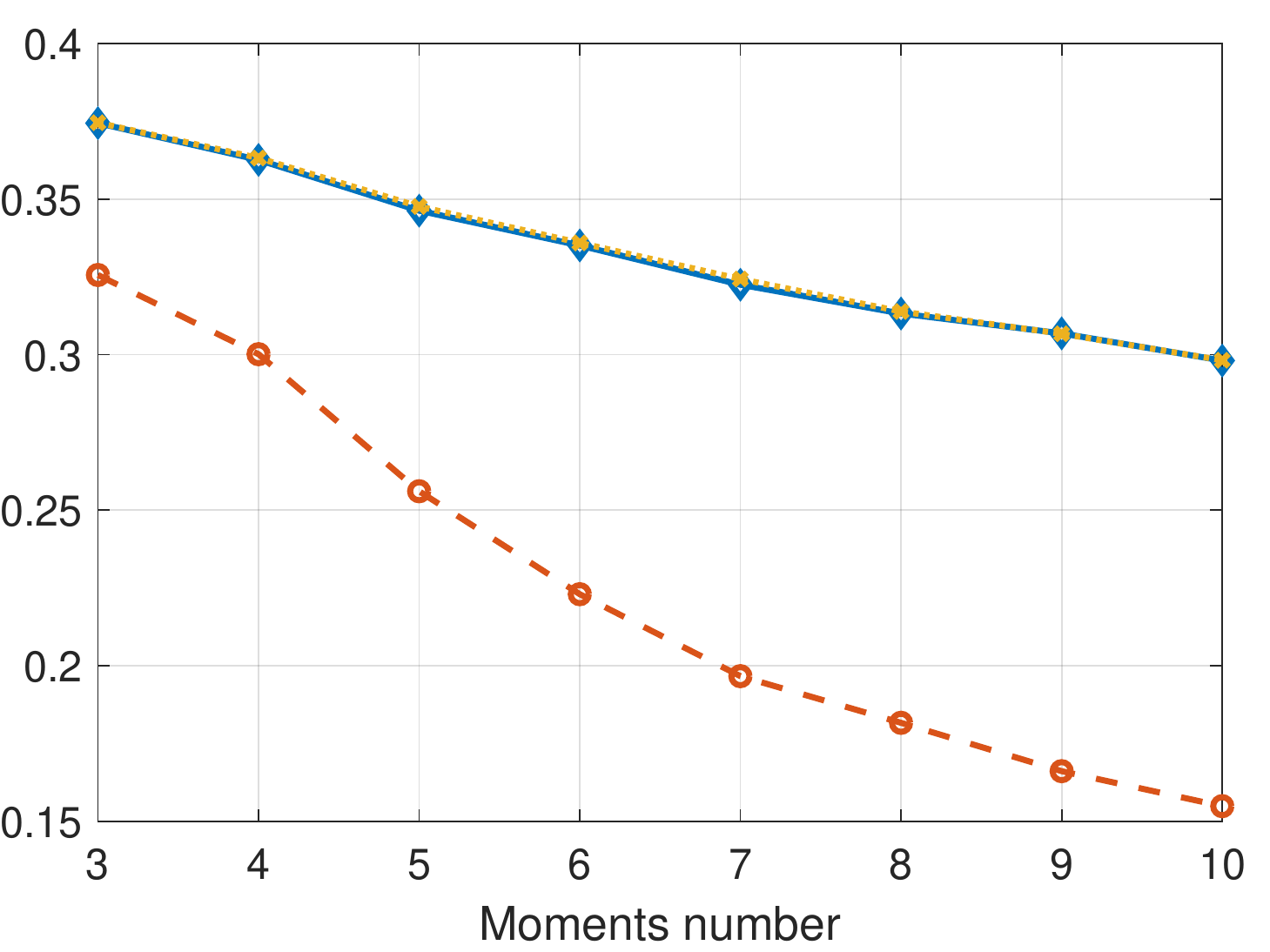}}
\bigskip
\centerline{\includegraphics[width=.48\textwidth]{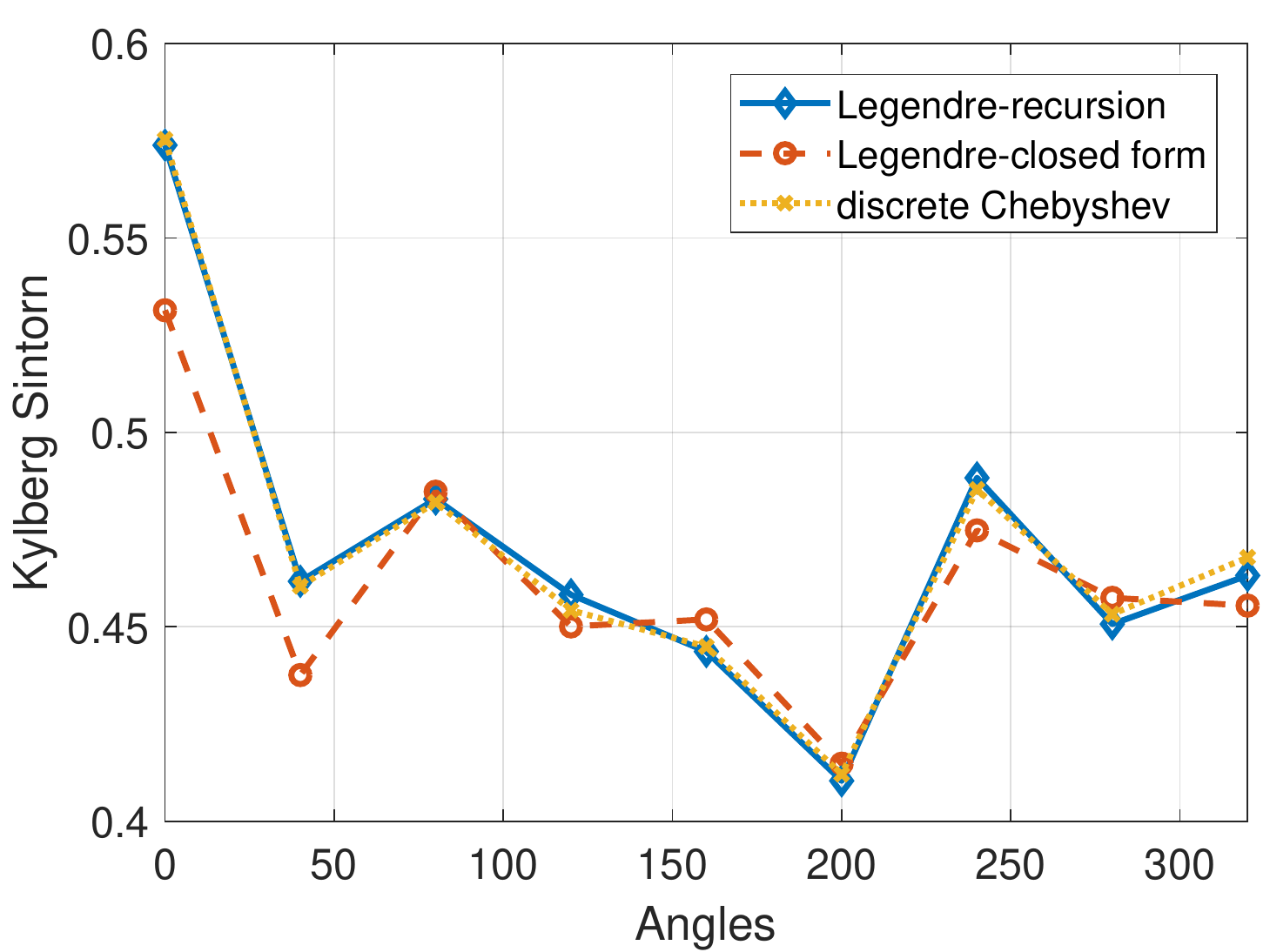} \hfill \includegraphics[width=.48\textwidth]{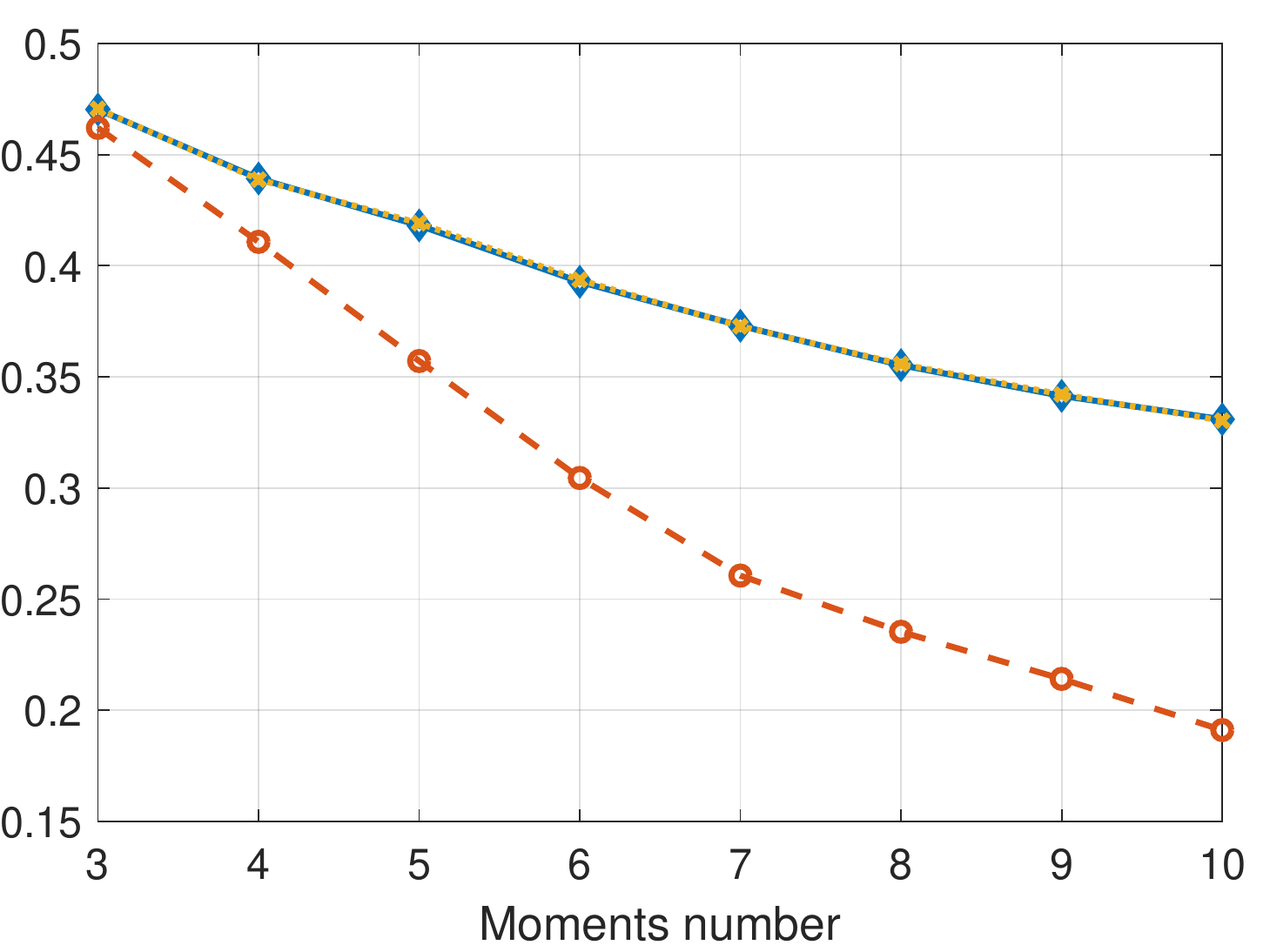}}
\bigskip
\centerline{\includegraphics[width=.48\textwidth]{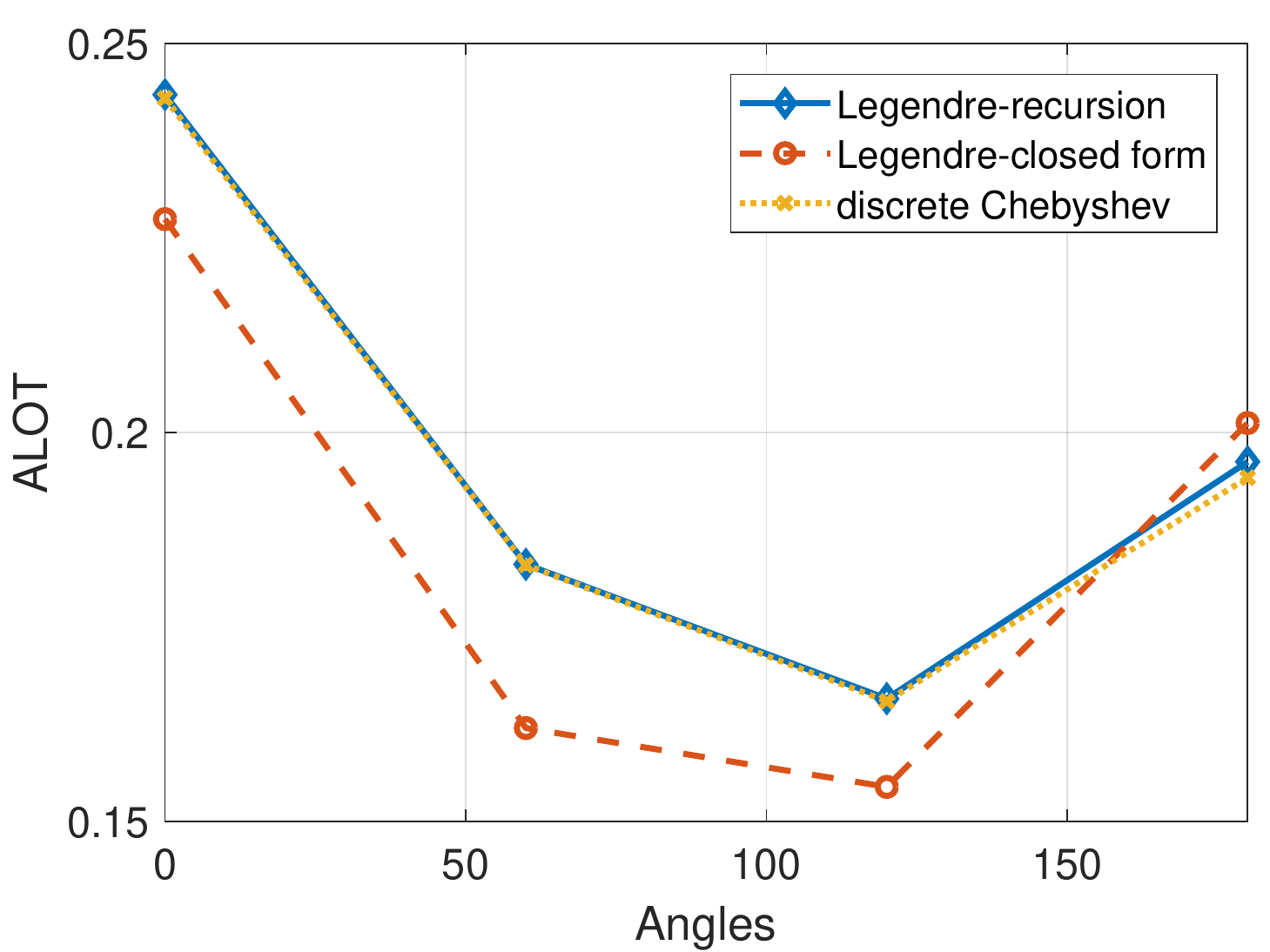} \hfill \includegraphics[width=.48\textwidth]{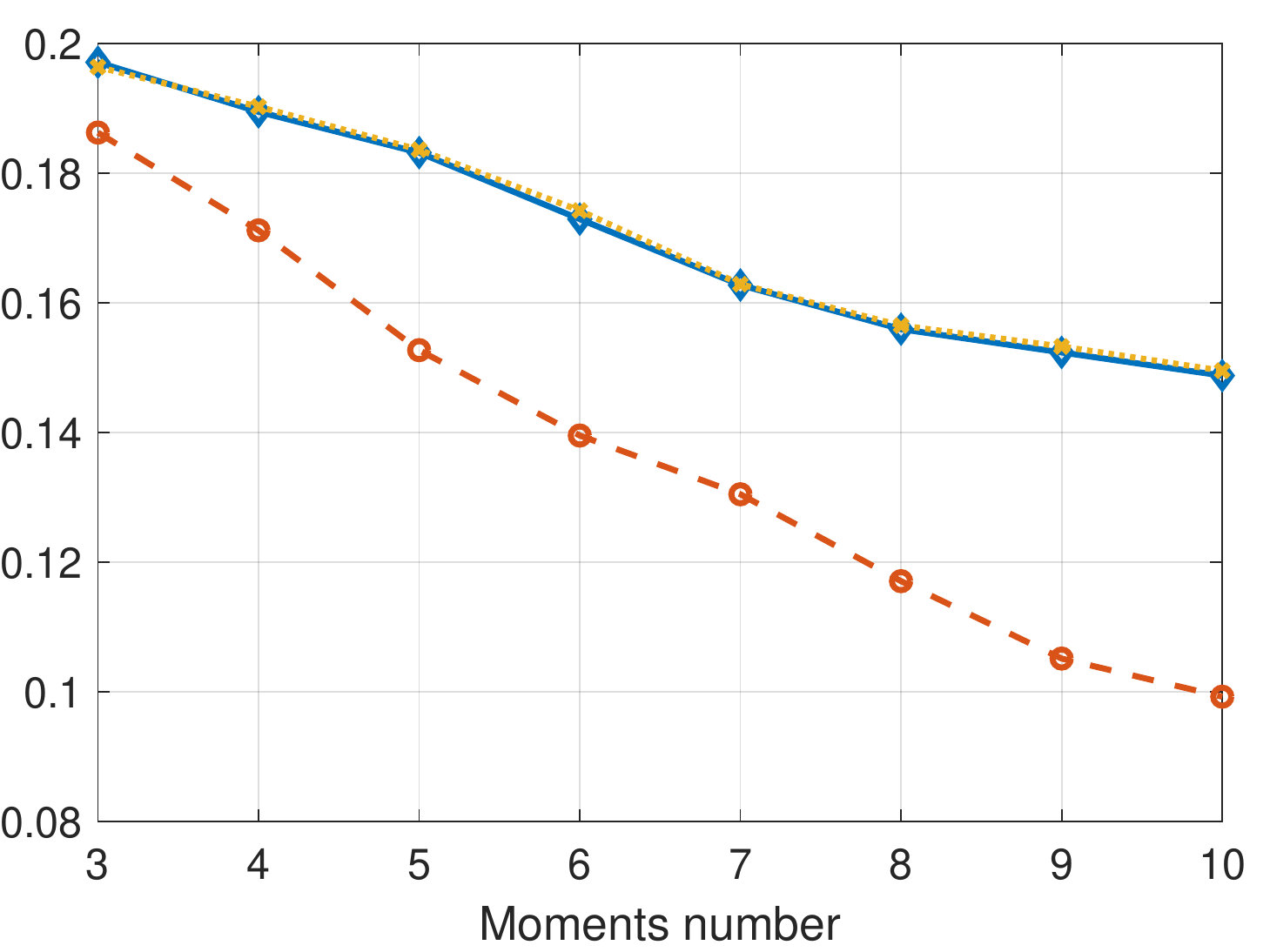}}
\caption{Vectorial, Kylberg Sintorn, and ALOT databases: accuracy index 
\eqref{accuracy} versus image orientation (left) and number of moments (right).}
\label{fig:performances2}
\end{center}
\end{figure}

\begin{figure}[!tb]
\begin{center}
\includegraphics[width=.488\textwidth]{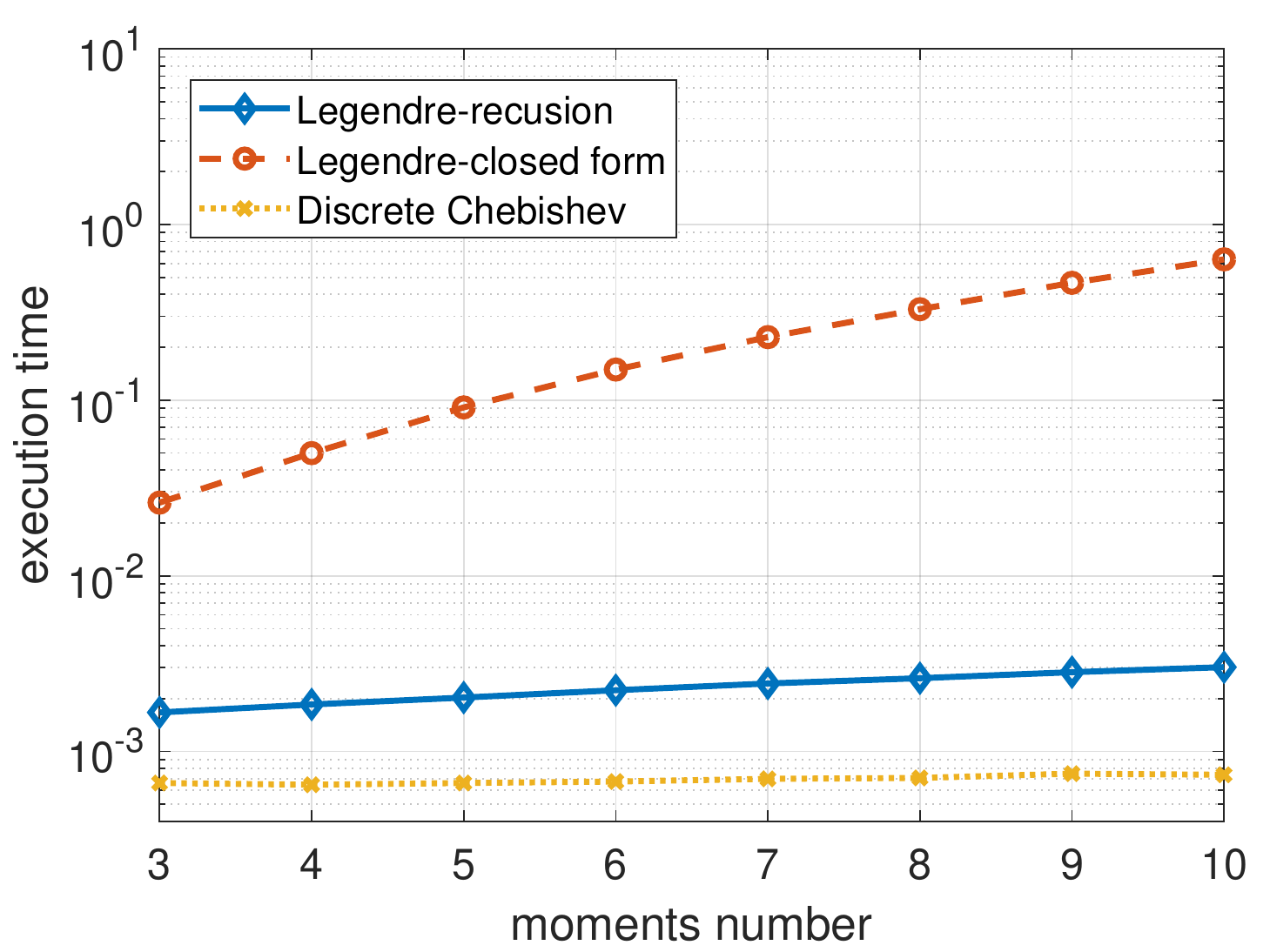}
\includegraphics[width=.47\textwidth]{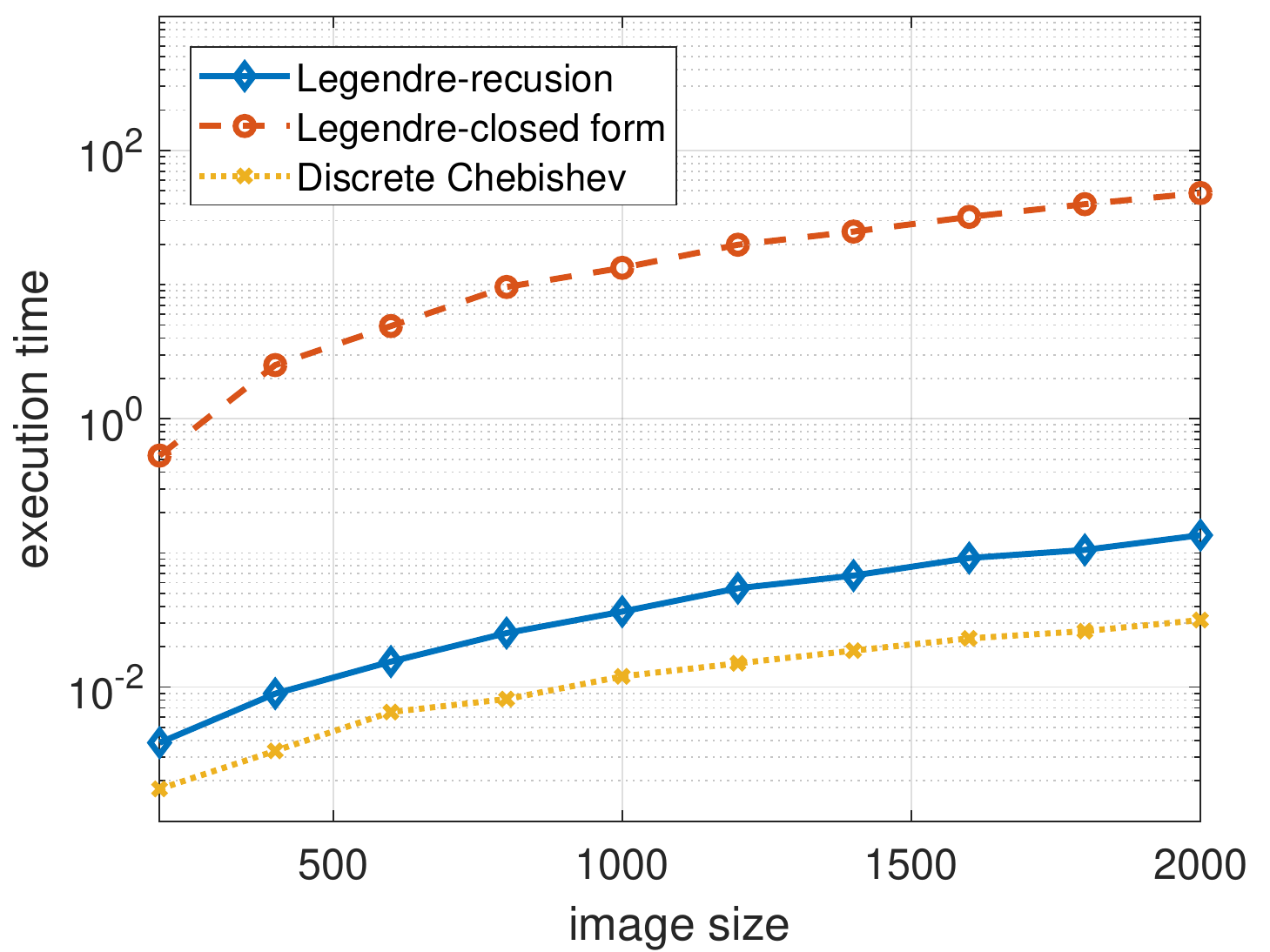}
\caption{Execution time when the number of moments takes the values
$n=3,4,\ldots,10$ (left plot), and with a fixed number of moments ($n=10$,
right plot) and image size varying from $200\times200$ up to $2000\times2000$.}
\label{fig:timing}
\end{center}
\end{figure}

The results of this experiment are reported in Figures~\ref{fig:performances}
and~\ref{fig:performances2}.
The graphs on the left display the accuracies for the different algorithms with
a fixed number of moments (q=3) and
varying the orientation of the image; the graphs on the right report the same
quantity versus the variation of the number of moments.
As it can be observed, the results obtained using the Legendre polynomials and
the discrete Chebyshev polynomials computed by recursion formulas 
present roughly the same accuracy values, and appear to be much more stable
than the closed form representation, for what concerns the orientation.
Indeed, the computational scheme used by the two approaches based on recurrence
relations leads to a better performance for the classification task.
The same trend can be observed also when the number of moments varies; see the
graphs on the right in the same figures.
In this case, we reported the average accuracy value resulting from all the
orientations.
The computing time is reported in Figure~\ref{fig:timing}, showing again the
superiority of the proposed approach in terms of computational complexity.

Figures~\ref{fig:performances} and~\ref{fig:performances2} also show that,
unlike the reconstruction phase, the performance in classification decreases
when the number of moments grows and, therefore, as the related descriptors
specialize.
Even though the recurrent formulation produces a better performance in
classification than the closed form, it is still insufficient for accurate
texture classification.
For this reason, we performed a further experiment using an approach
already proposed in \cite{DiRuberto2017}, where
the moments were computed starting from a different
representation of the images. Indeed, higher level features were created using
the Gray Level Co-occurrence Matrices (GLCM) \cite{Benco,Chen,Haralick,Mitrea}
and the Local Binary Pattern (LBP) \cite{Ojala96}, since it was observed that
using these approaches (in particular, the first one) the invariant moments are
more discriminative.

\begin{figure}[!t]
\begin{center}
\centerline{\includegraphics[width=.48\textwidth]{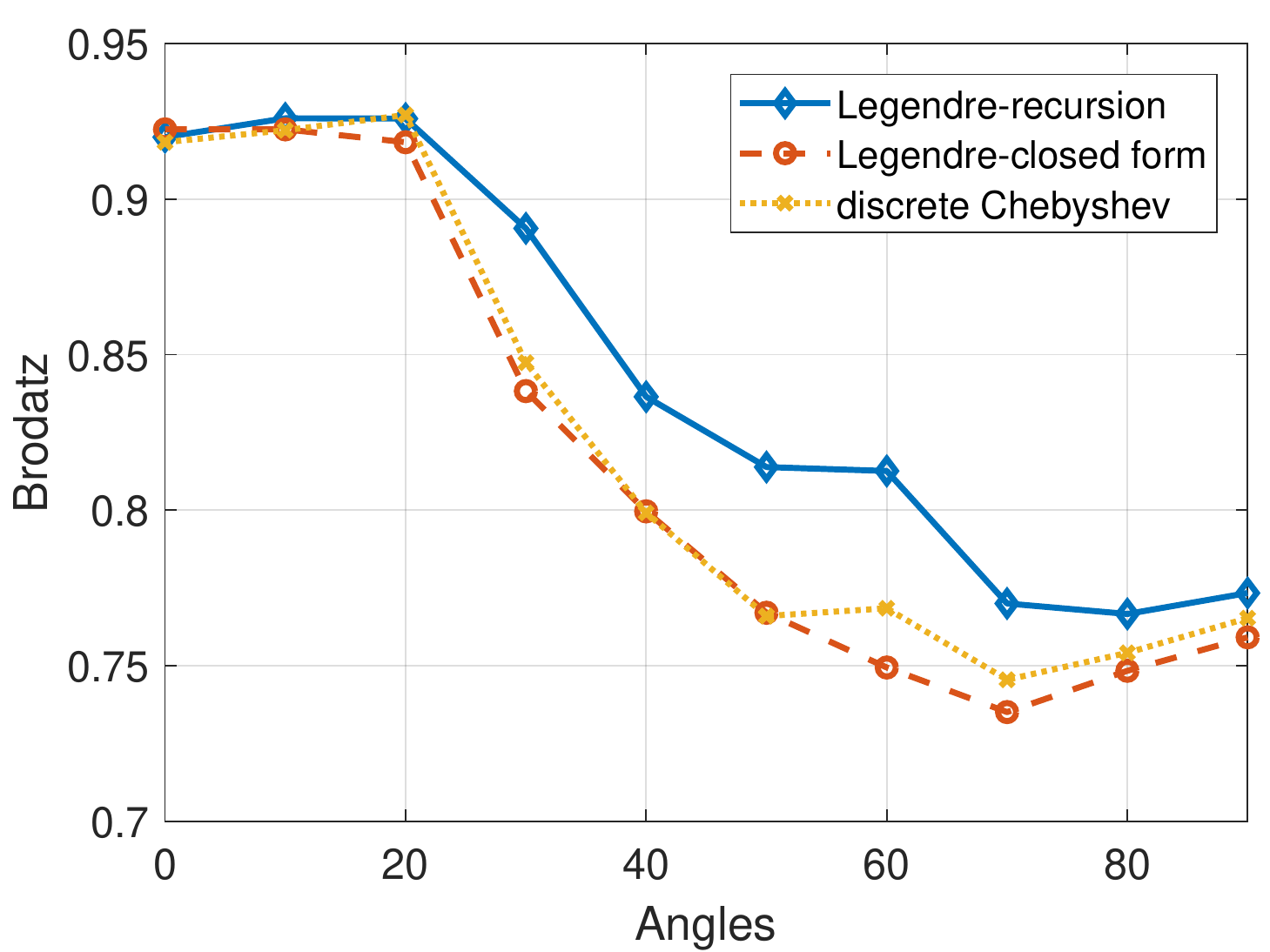} \hfill \includegraphics[width=.48\textwidth]{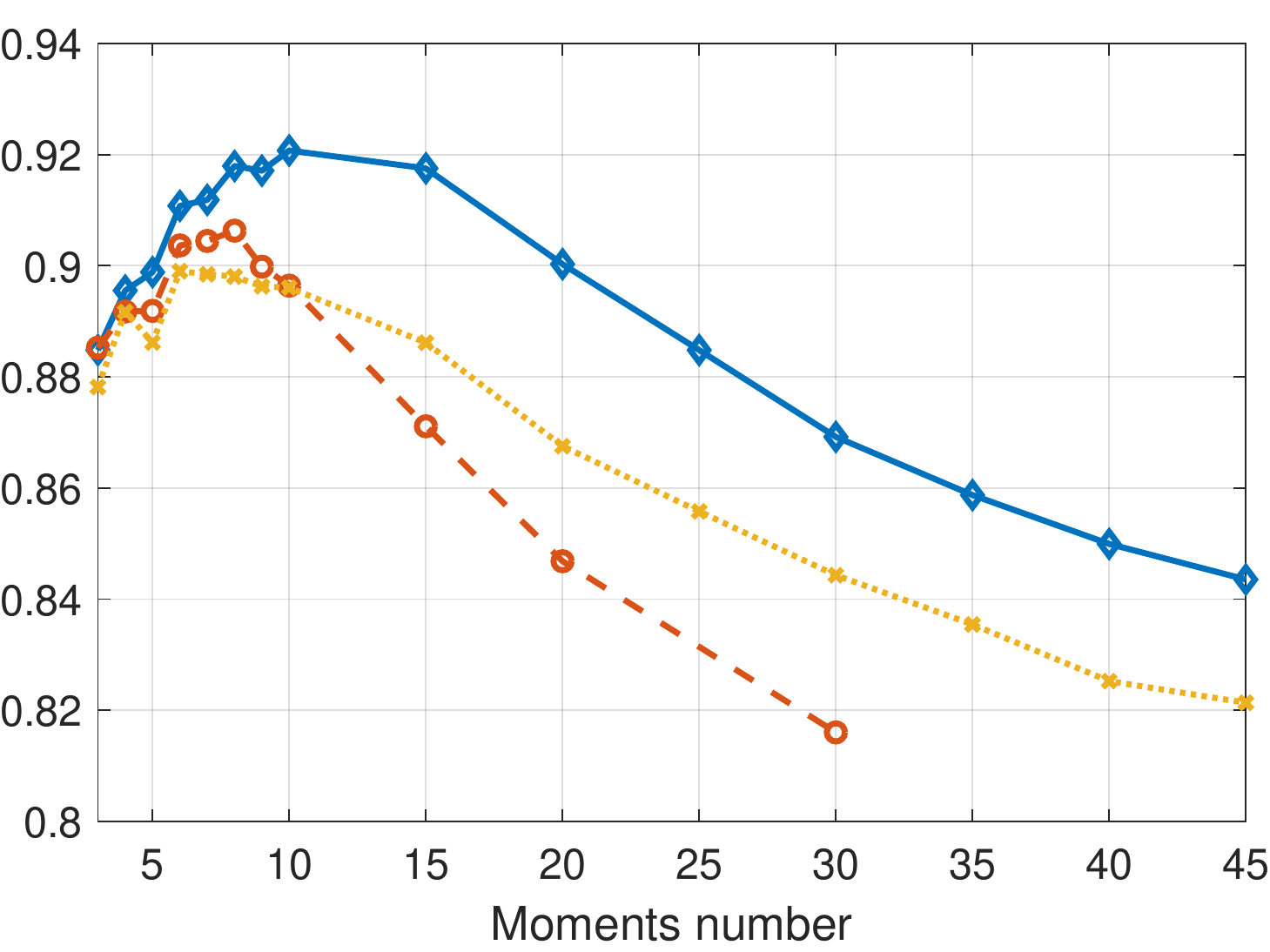}}
\bigskip
\centerline{\includegraphics[width=.48\textwidth]{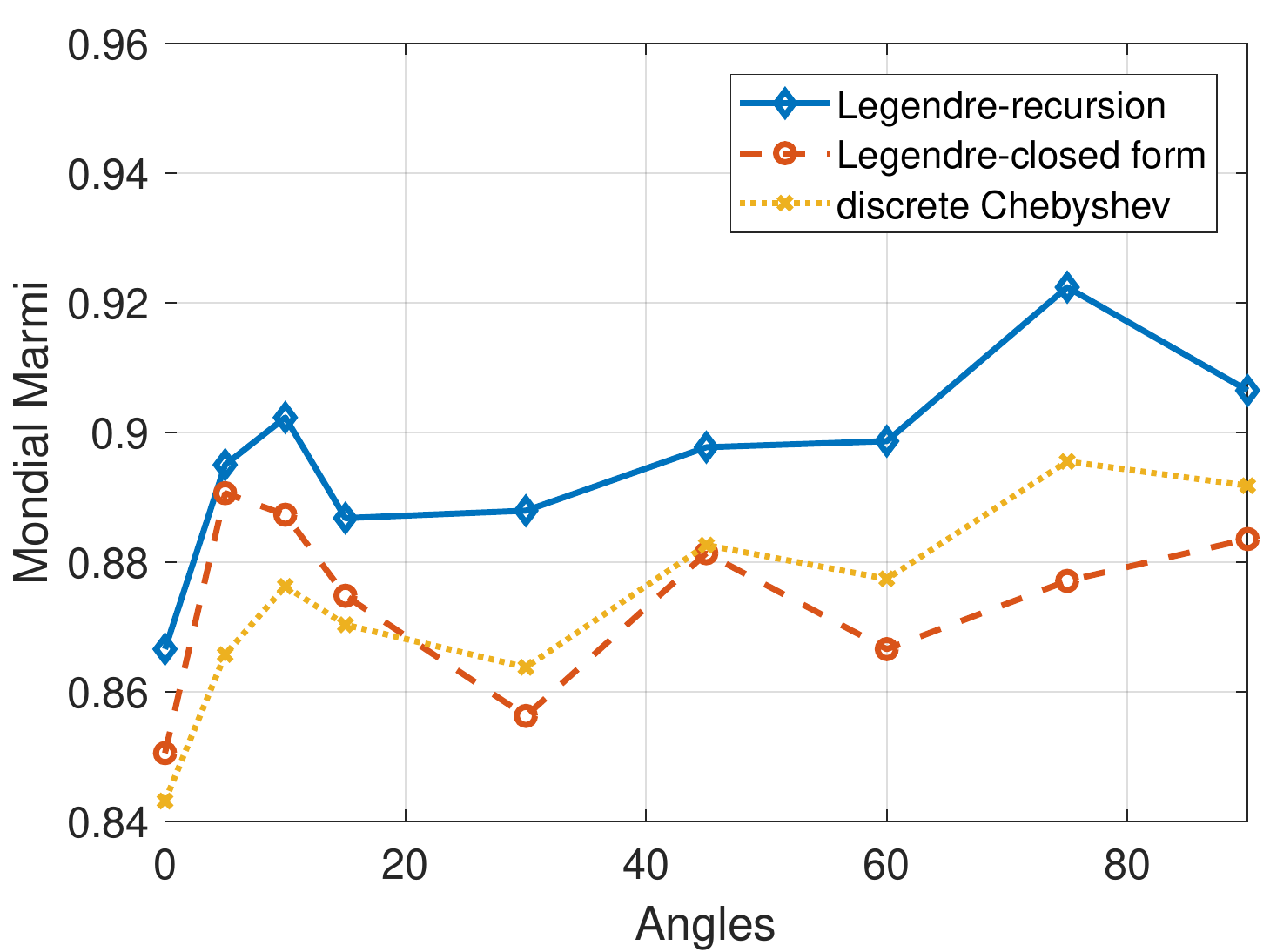} \hfill \includegraphics[width=.48\textwidth]{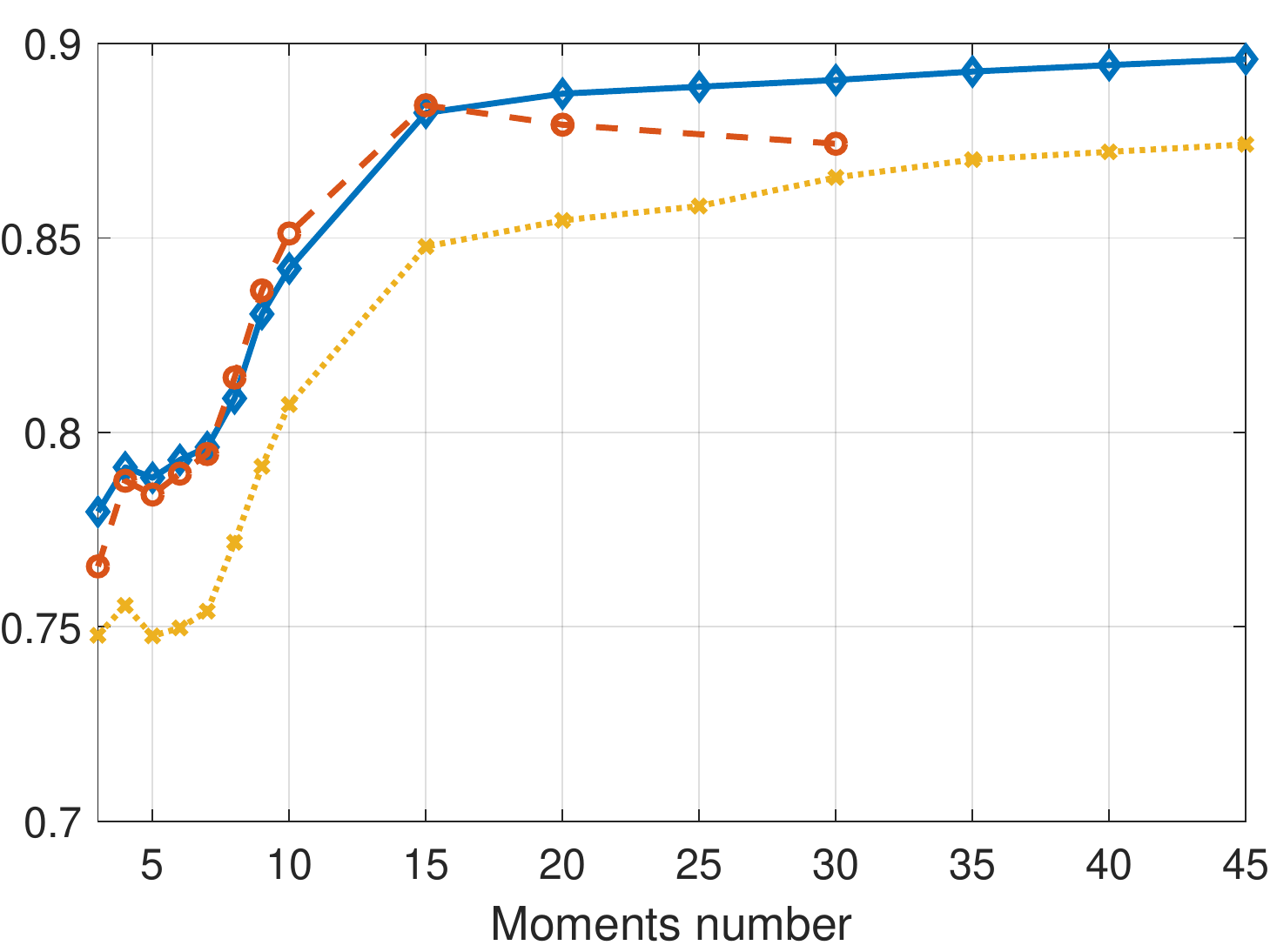}}
\bigskip
\centerline{\includegraphics[width=.48\textwidth]{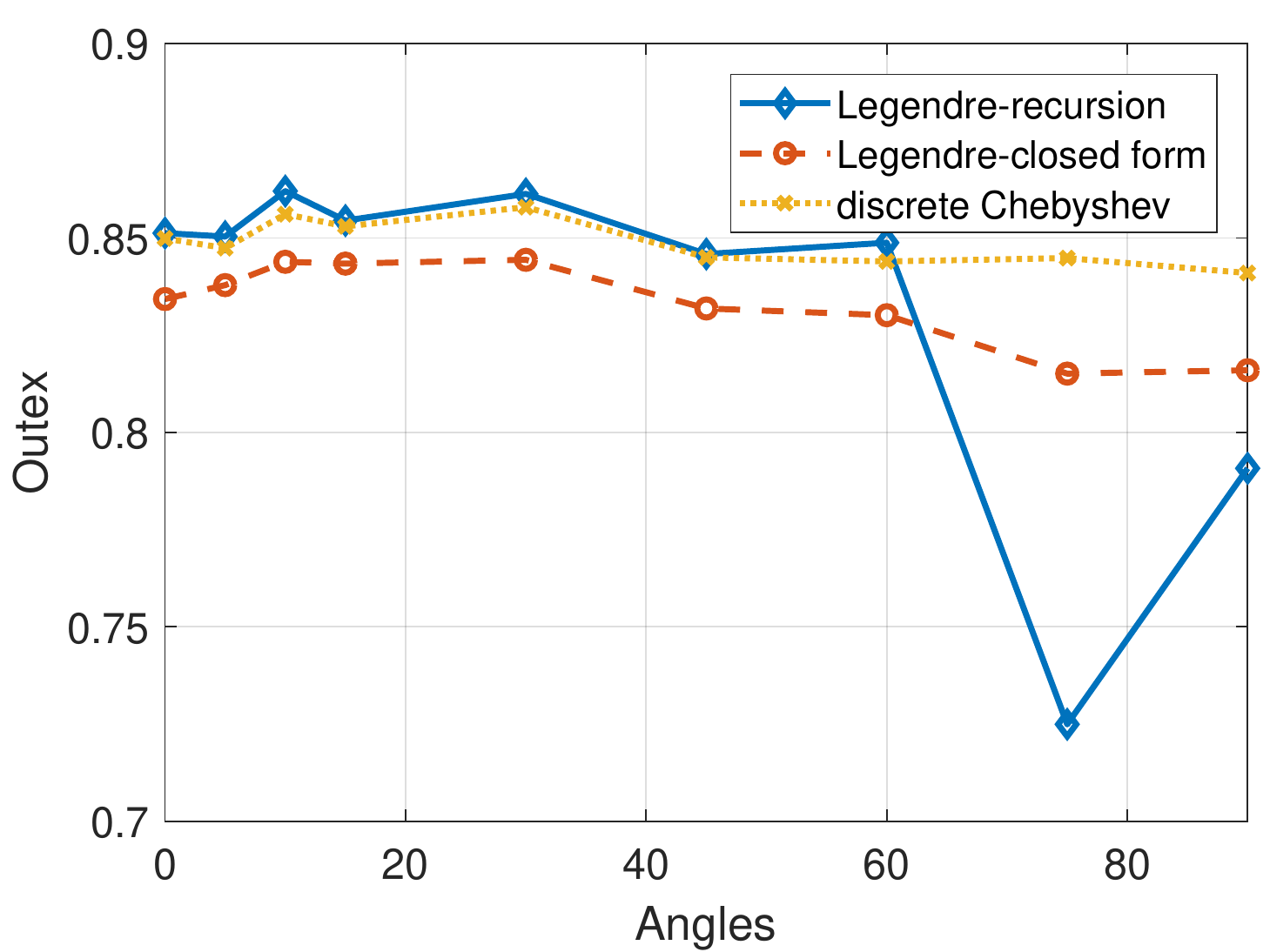} \hfill \includegraphics[width=.48\textwidth]{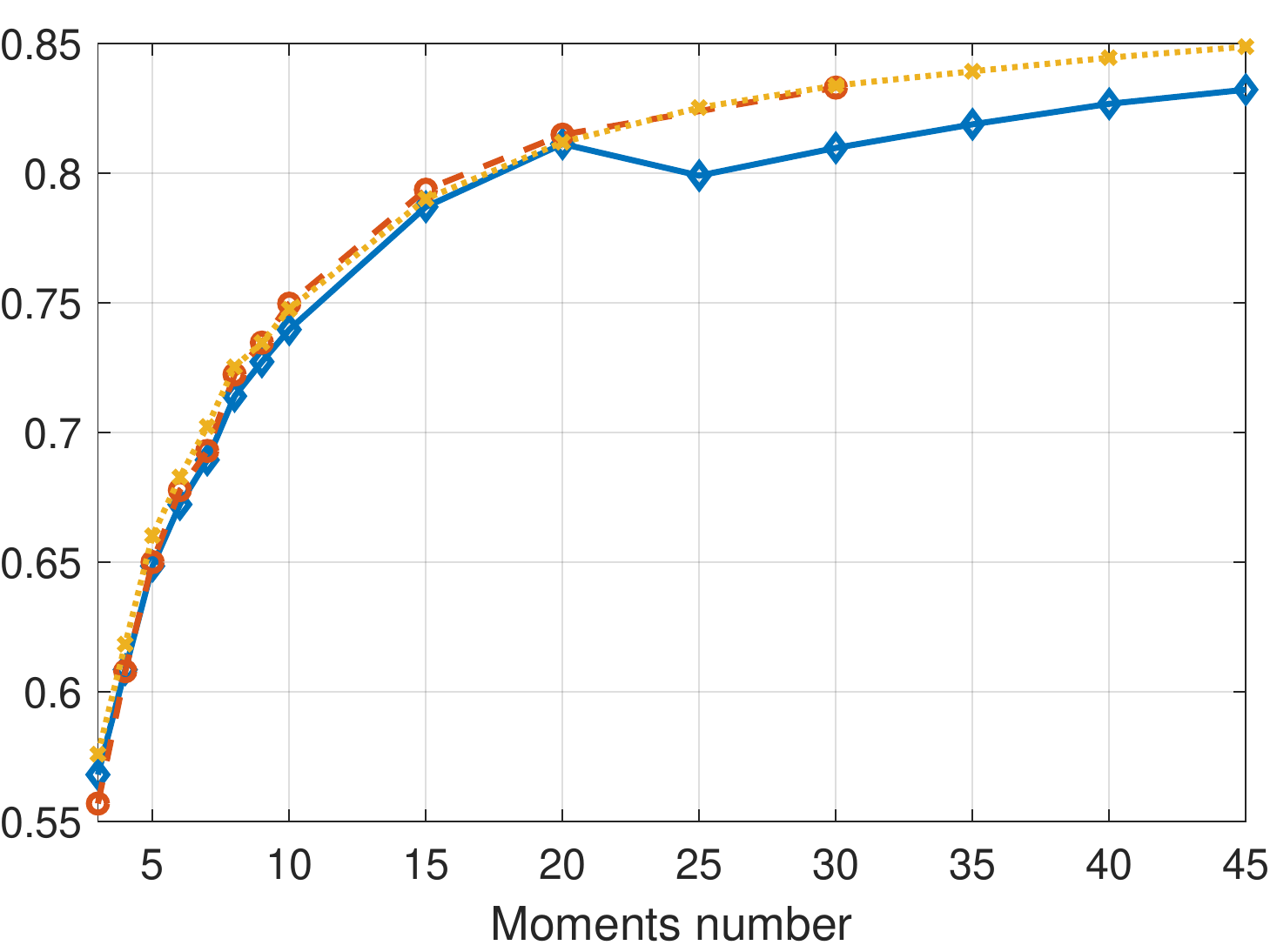}}
\caption{Brodatz, Mondial Marmi, and Outex databases: accuracy index
\eqref{accuracy} versus image orientation (left) and number of moments (right),
when the moments are extracted from the GLCM.}
\label{fig:performancesGLCM}
\end{center}
\end{figure}

\begin{figure}[!t]
\begin{center}
\centerline{\includegraphics[width=.48\textwidth]{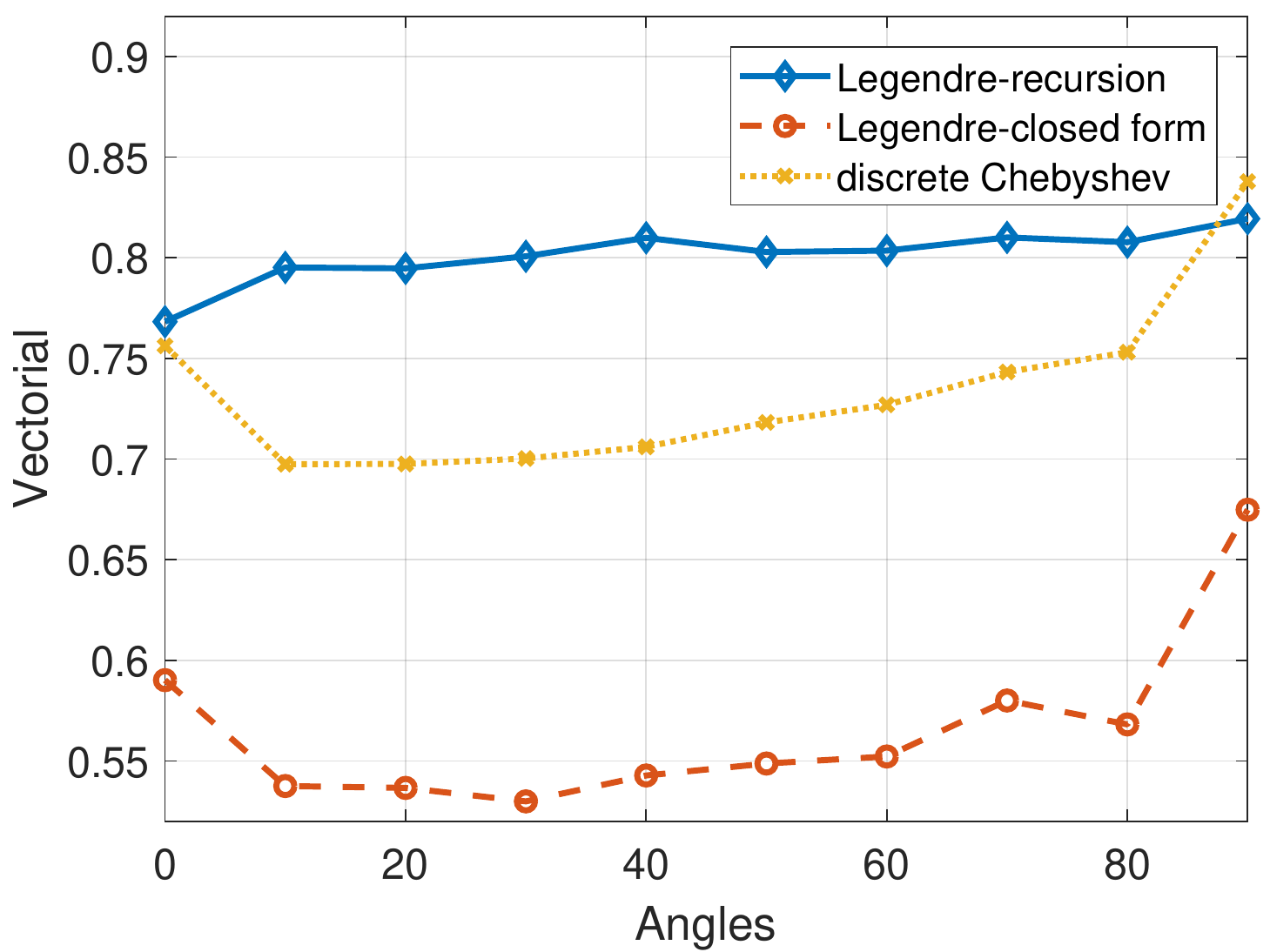} \hfill \includegraphics[width=.48\textwidth]{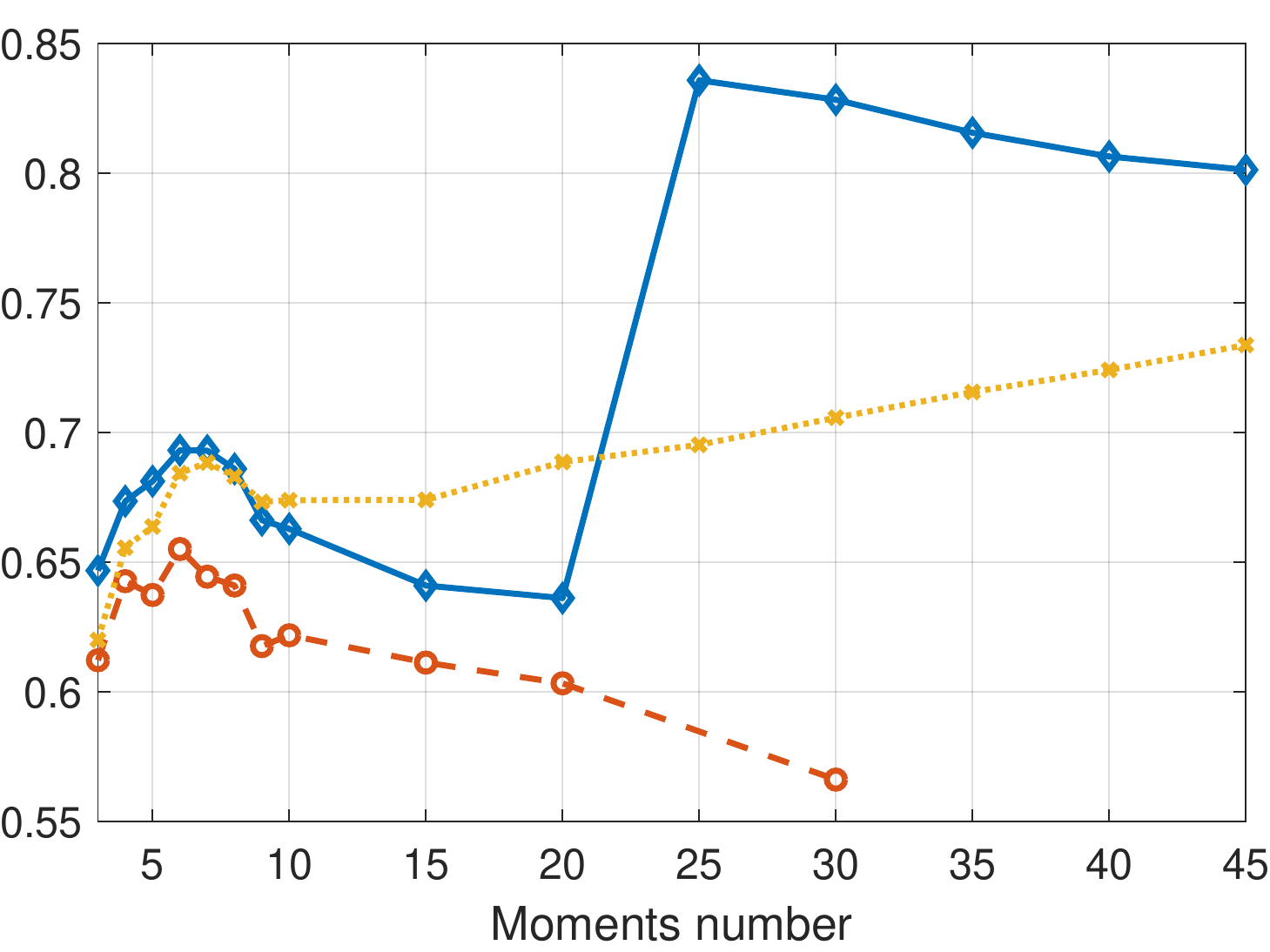}}
\bigskip
\centerline{\includegraphics[width=.48\textwidth]{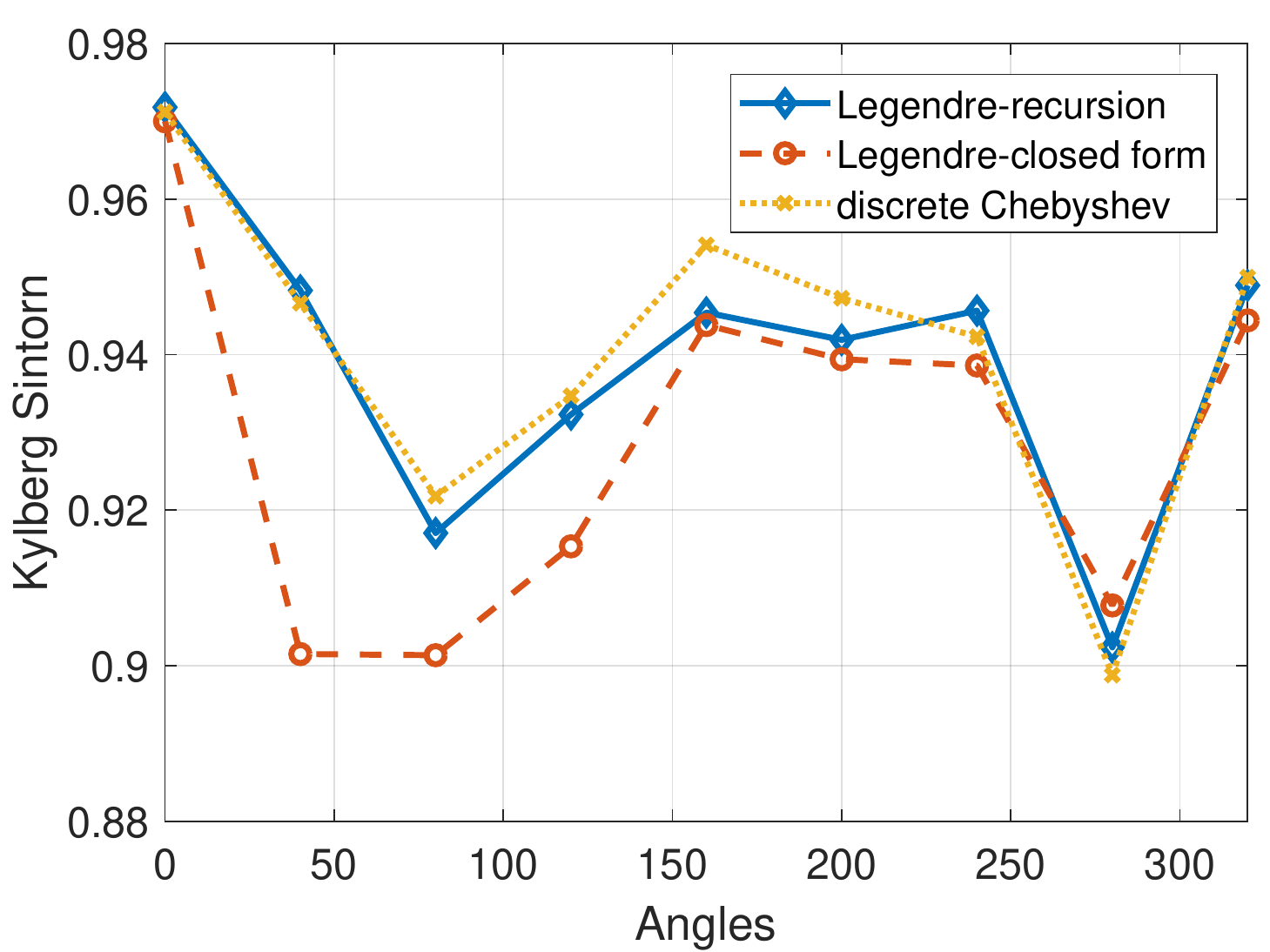} \hfill \includegraphics[width=.48\textwidth]{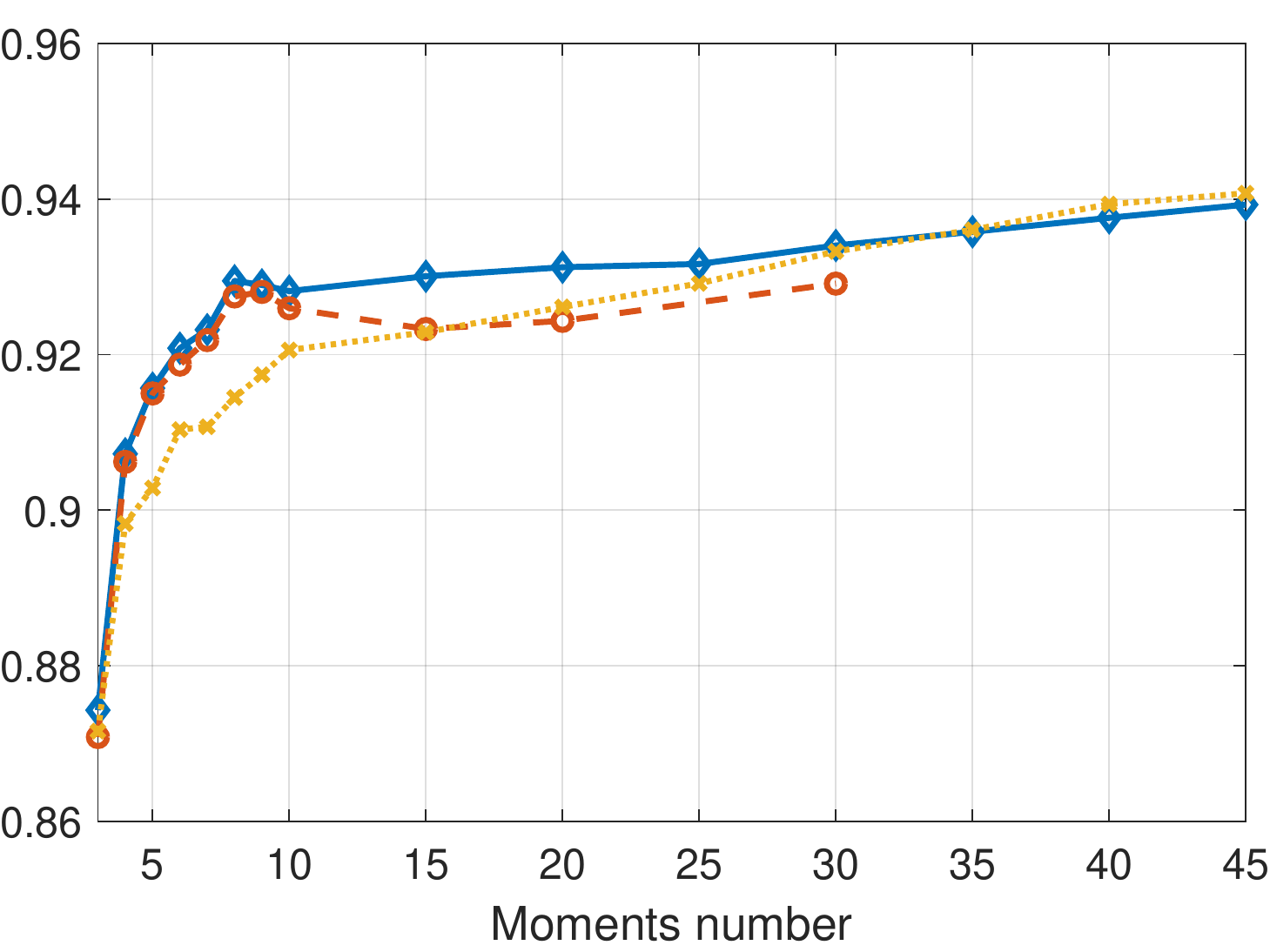}}
\bigskip
\centerline{\includegraphics[width=.48\textwidth]{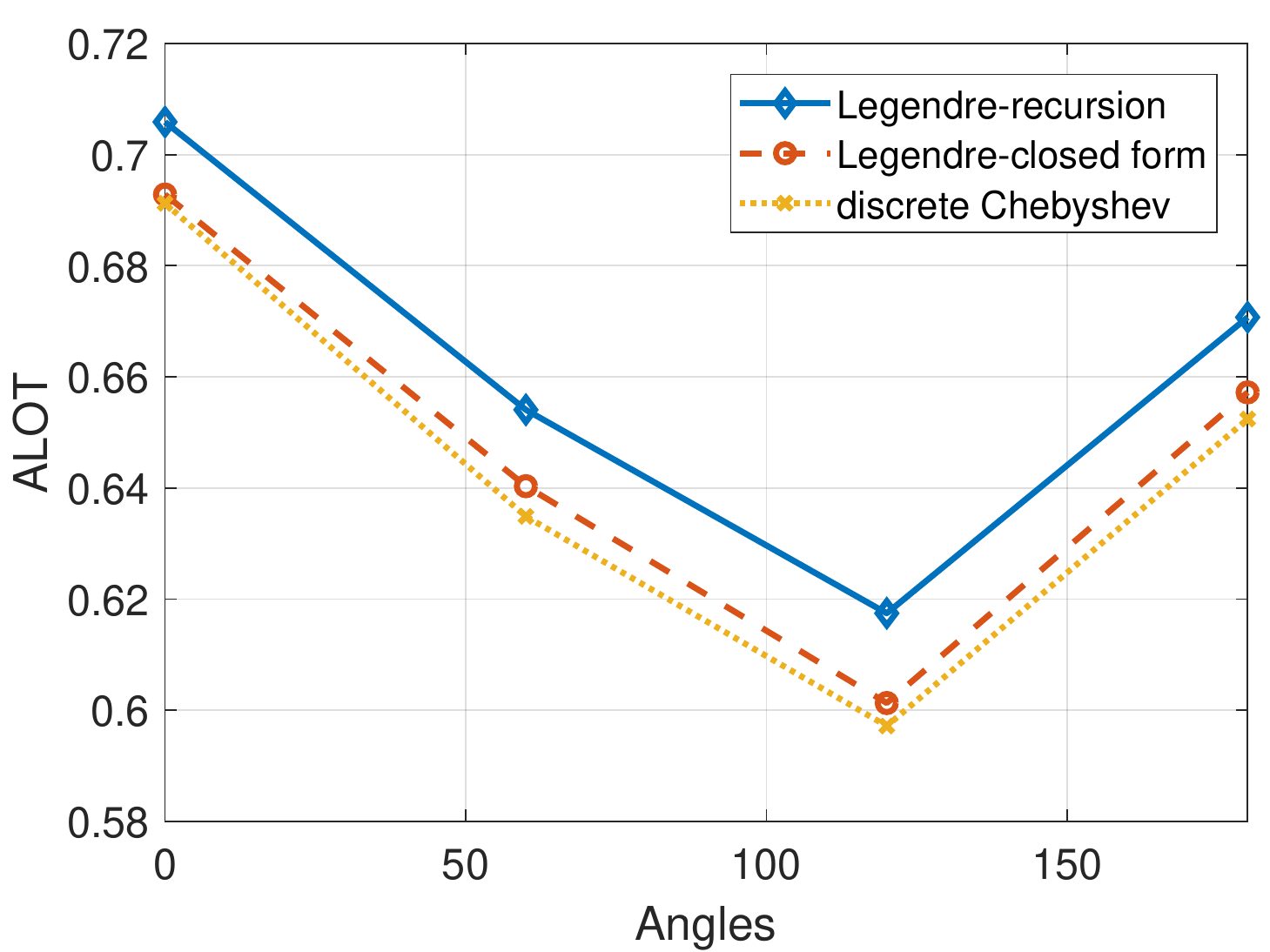} \hfill \includegraphics[width=.48\textwidth]{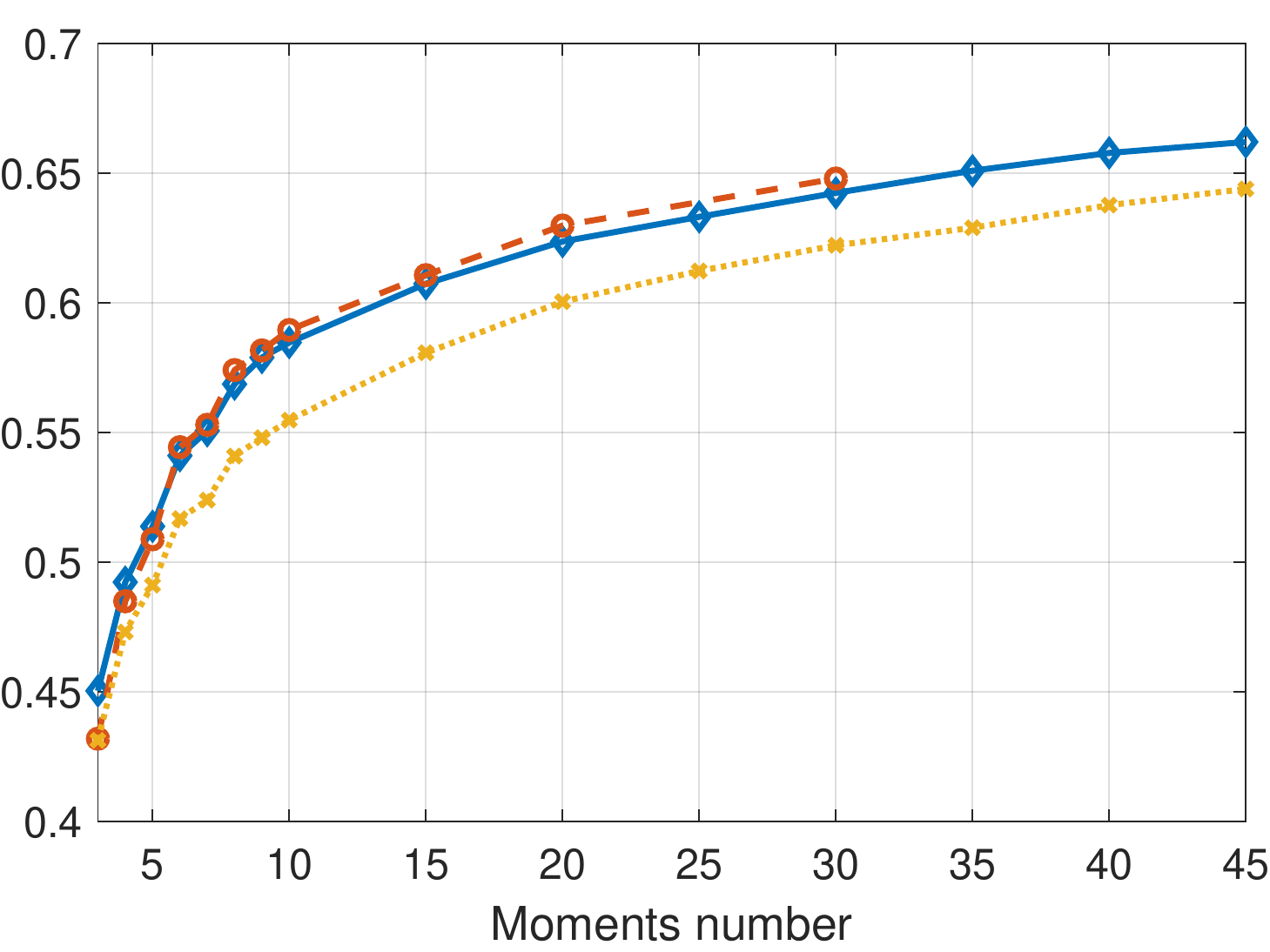}}
\caption{Vectorial, Kylberg Sintorn, and ALOT databases: accuracy index 
\eqref{accuracy} versus image orientation (left) and number of moments (right),
when the moments are extracted from the GLCM.}
\label{fig:performancesGLCM2}
\end{center}
\end{figure}

Thus, we computed the moments starting from the GLCM generated with a distance
value $d=1$ and angles $0^\circ,45^\circ,90^\circ$, and $135^\circ$. The final
feature vector was created by concatenating the resulting moments. Since the
GLCM has been computed from 4 angles, the final feature vector is 4 time
larger than the previous one. Indeed, now the feature vector size ranges from
40 to 4324, with a number of moments ranging from 3 to 45, while the size of
the feature vector obtained by extracting the moments directly from the images
ranges from 10, when the number of moments is 3, to 66, when the number of
moments is 10. The results of this comparison are reported in
Figures~\ref{fig:performancesGLCM} and~\ref{fig:performancesGLCM2}. 

The execution time for the whole process of feature extraction is obviously
larger, since it also includes the time needed to compute the GLCMs. Anyway,
the execution time just increases by a constant value, since the GLCMs have a
fixed size, typically $256\times256$.
Thus, the trend remains essentially the same as depicted in
Figure~\ref{fig:timing}. At the same time, the performance in classification
greatly improves and the classification accuracy increases with the
number of moments. This can be observed in Figures~\ref{fig:performancesGLCM}
and~\ref{fig:performancesGLCM2}, where we compare the effects on accuracy of 
image rotation with a fixed number of moments (q=45) (graphs on the left) and
of the number of moments (on the right).

\begin{figure}[!t]
\begin{center}
\centerline{\includegraphics[width=.48\textwidth]{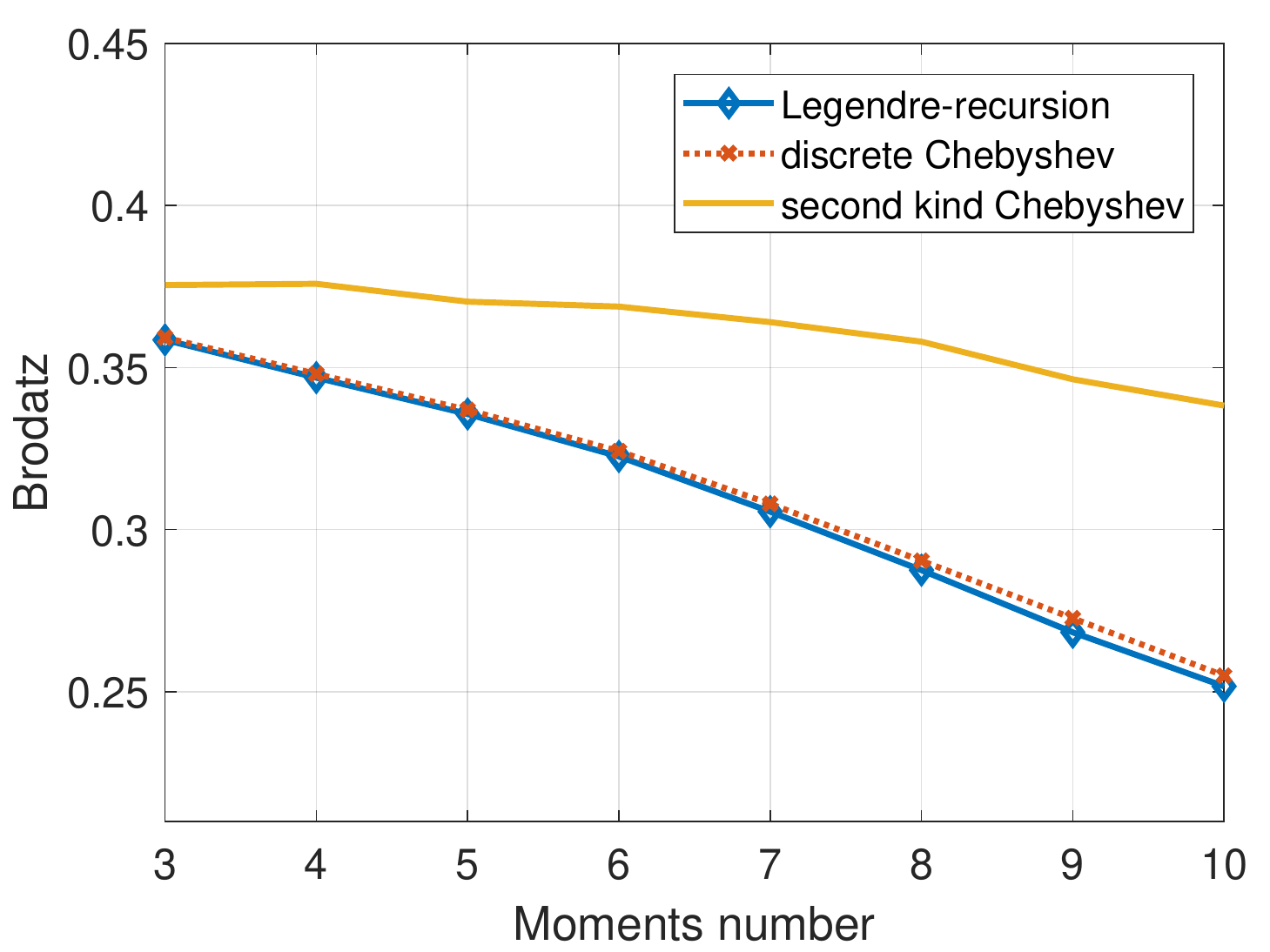} \hfill \includegraphics[width=.48\textwidth]{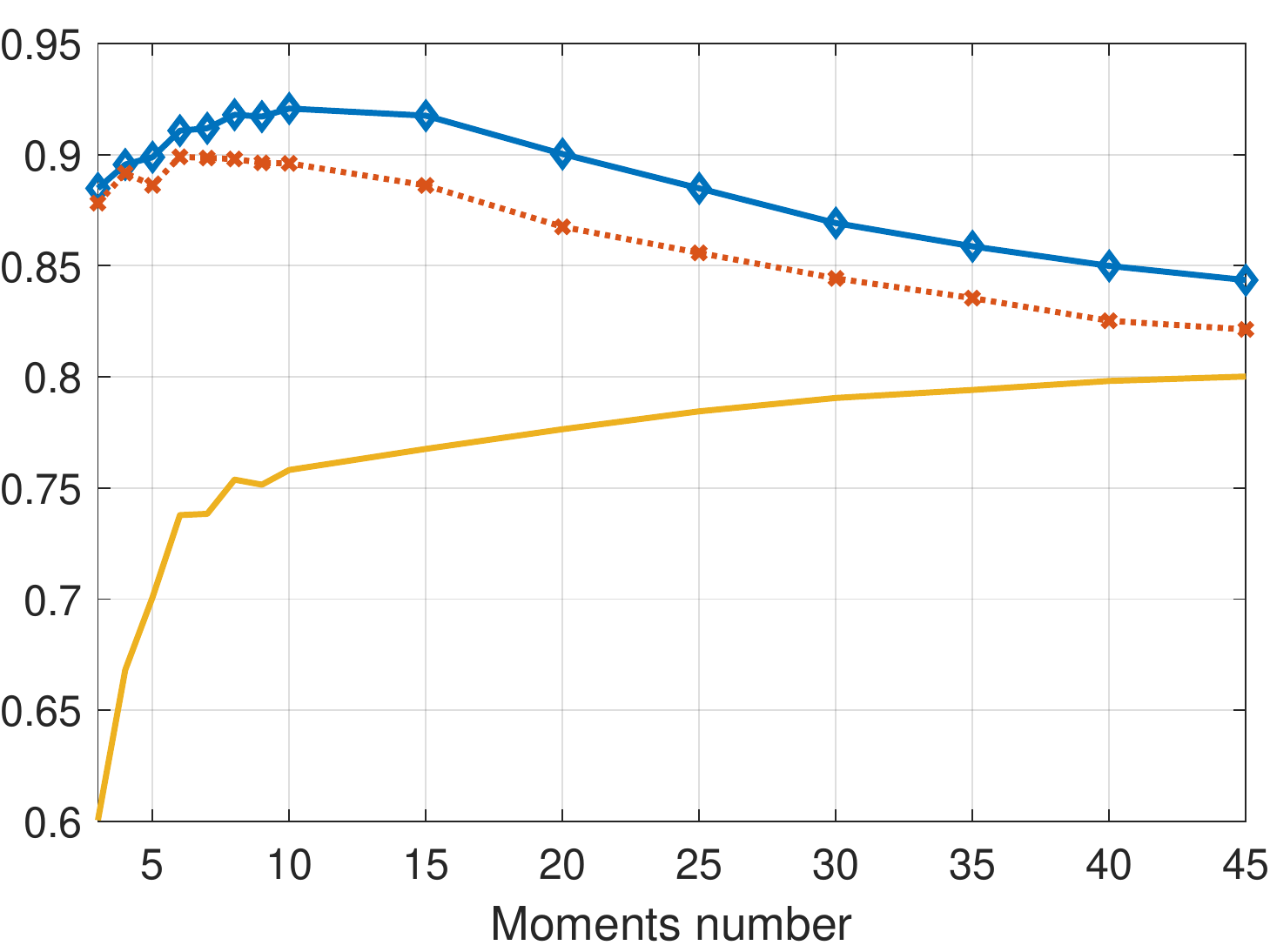}}
\bigskip
\centerline{\includegraphics[width=.48\textwidth]{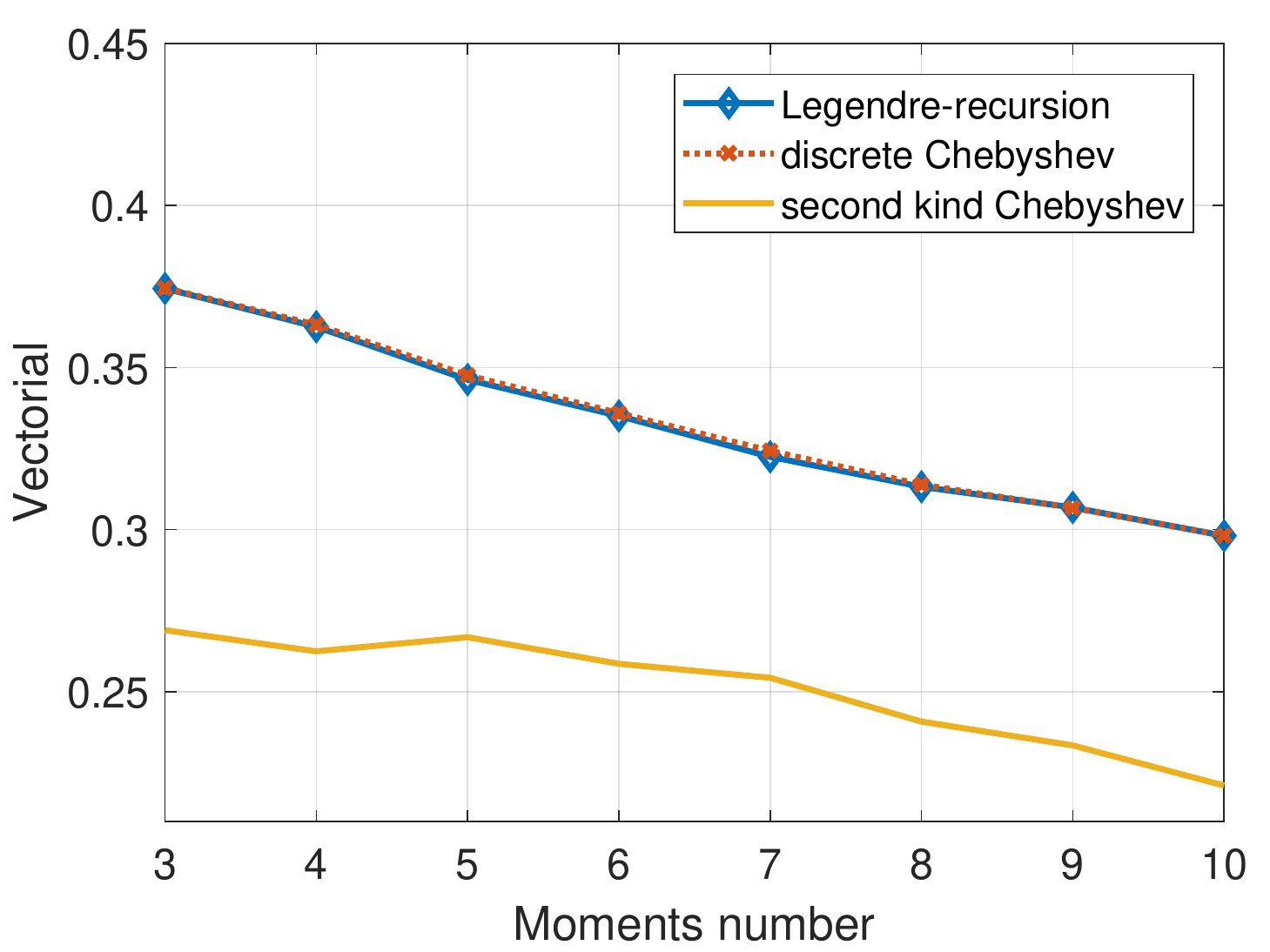} \hfill \includegraphics[width=.48\textwidth]{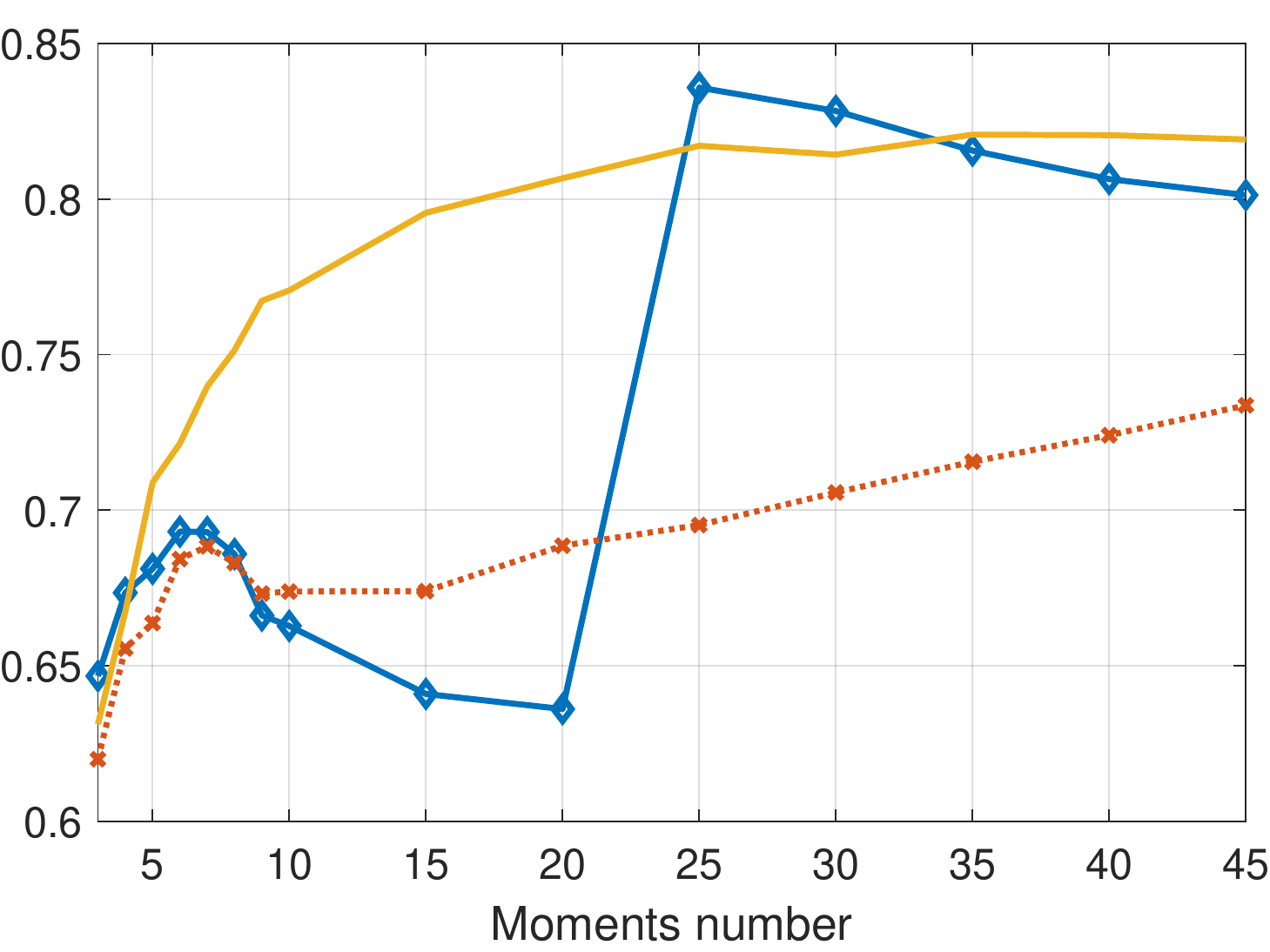}}
\caption{Brodatz and Vectorial databases: comparison between second kind
Chebyshev, Legendre, and discrete Chebyshev polynomials. The accuracy index
\eqref{accuracy} is plotted against the number of moments. The graphs in the
left column concern moments extracted from the image, while in the right
column graphs the moments are extracted from the GLCM.}
\label{fig:performancesCheb2}
\end{center}
\end{figure}

Second kind Chebyshev (Cheb2) polynomials \eqref{cheby2} are orthogonal with
respect to the inner product \eqref{innprod}, with $w(x)=\sqrt{1-x^2}$.
This weight function achieves its maximum at the center of $[-1,1]$, while
it drops to zero at the endpoints. This means that the computation of the
moments \eqref{Lmom_cont} emphasizes the central part of the image and reduces
the importance of its border area. This effect may be useful when classifying
images having a particular structure. This fact is outlined in
Figure~\ref{fig:performancesCheb2}, where we display the accuracy index
\eqref{accuracy} versus the number of moments used for classification, for two
datasets. The results displayed in the left column are obtained by extracting
the moments directly from the images. In this case, the accuracy index for the
Brodatz database improves when Cheb2 polynomials are employed, while their
performance is much worse than Legendre and discrete Chebyshev polynomials for
the Vectorial database. On the contrary, when the moments are extracted from
GLCMs, Cheb2 polynomials overall performance is the best for Vectorial, while
it is not very good for Brodatz. We believe that weighted orthogonal moments
may have a great potential in image classification, but a much deeper analysis
is needed in order to understand when they are preferable to unweighted
moments. This will be the subject of future work.

To further highlight the performance of the approaches studied in this paper,
we compared them to some widely used state-of-the-art descriptors for texture
classification. We computed the rotation invariant GLCM features as proposed in
\cite{DiRuberto2017b}, the rotation invariant LBP (LBP-RI) \cite{Ojala02} with
a $3\times3$ pixels neighborhood and distance 1, and the Convolutional Neural
Network (CNN) features from three different well known network architectures.
The first architecture is AlexNet \cite{Krizhevsky12} that gained popularity
for its good performance in many classification tasks. It consists of 8 layers,
and we extracted the features from the activations of the second last fully
connected layer that produces a feature vector of size 4096. The second
architecture is the Vgg19Net \cite{Simonyan:2014aa}, which is quite similar to
AlexNet except for the number of layers being 19. Also in this case, we
extracted the features from the activations of the second
last fully connected layer that produces a feature vector of size 4096. The
third architecture, GoogLeNet \cite{Szegedy2015}, is quite different.
It presents a stack of 22 layers but, since it uses inception layers, the total
number of layers rises to 100. We extracted the features from the activations
of the last fully connected layer (the only one that can be used for feature
extraction) that leads to a feature vector of size 1000.

As in the previous comparison, we evaluated the effect of image rotation by
training the k-NN
classifier with features extracted from images acquired at orientation
$0^\circ$, and then we tested with features extracted from images acquired at
other orientations. To better illustrate the overall performance, we computed
the average results for all the angles. The outcome of this experiment is
reported in Table~\ref{tab:performancesALL}, where we compare the accuracy of
the algorithms investigated in this paper to the state-of-the-art approaches.
The best results for each dataset are displayed in bold face. 

\begin{table}[!t]
\caption{Comparison of classical invariant moments, and invariant moments
extracted from GLCM, to three state-of-the-art approaches: rotation invariant
GLCM Features, rotation invariant LBP, and CNN Features from AlexNet, Vgg19Net,
and GoogLeNet.}
\centering\tabcolsep=1mm
\begin{tabular}{lcccccc}
\toprule
Texture descriptor & \multicolumn{6}{c}{Dataset} \\ 
& Brodatz & Mondial & Outex & Vectorial & Kylberg & ALOT \\ 
\midrule
Legendre-closed form & 32.8  & 48.3  & 26.0  & 32.6  & 46.2  & 18.6 \\ 
Legendre-recursion & 35.9  & 54.2  & 25.3  & 37.4  & 47.0  & 19.7 \\ 
Discrete Chebyshev & 35.9  & 54.1  & 25.3  & 37.5  & 47.1  & 19.6 \\ 
Legendre-c.f. GLCM & 90.6  & 88.4  & 83.3  & 65.5  & 92.9  & 64.8 \\ 
Legendre-rec. GLCM & \textbf{92.1}  & \textbf{89.6}  & 83.2  & 83.6  & 93.9  & 66.2 \\ 
Discr. Cheb. GLCM & 89.9  & 87.4  & \textbf{84.9}  & 73.4  & 94.1  & 64.4 \\ 
GLCM Features & 66.2  & 65.6  & 61.9  & 50.5  & 85.7  & 40.7 \\ 
LBP-RI & 83.3  & 82.5  & 74.6  & 70.2  & 89.9  & 70.6 \\ 
AlexNet Features & 85.7  & 80.5  & 78.5  &90.5  & 91.8  & 80.7 \\
Vgg19Net Features & 86.9  & 76.7  & 81.3  &   \textbf{96.3}  & \textbf{95.1}  & \textbf{91.0} \\
GoogLeNet Features & 87.6  & 70.8  & 72.2  & 95.0  & 94.0  & 82.6 \\ 
\bottomrule
\end{tabular}\label{tab:performancesALL}
\end{table}

As it can be observed, the proposed formulation outperforms the
state-of-the-art approaches in 3 datasets out of 6 (Brodatz, Mondial, and
Outex), and produces essentially equivalent results for the Kylberg Sintorn
dataset.
On the contrary, the CNN features descriptors perform sensibly better on the
Vectorial and ALOT datasets.
This is probably due to the fact that these datasets present a very large
number of classes and, since the CNN features are extracted from networks that
learn from millions of images, in this case they are more effective in
describing the image content.

There is another possible interpretation for the poor performance of the
orthogonal moments on the Vectorial and ALOT datasets.
These two collections of images, as it can be observed in
Figure~\ref{fig:textures}, exhibit rather regular shapes and smooth textures.
The other datasets, instead, contain fine textures, characterized by the
presence of high frequencies.
It is well known that when a function is approximated in the least squares
sense by a sum of orthogonal functions, the decay speed of the expansion
coefficients increases in correspondence to an increase in the function
regularity.
So, there is the possibility that the decay in the moments extracted from the
Vectorial and ALOT datasets makes it difficult to perform an effective moment
matching during classification.
This fact needs a deeper investigation and understanding, which may lead to 
an effective technique to improve the classification performance of orthogonal
moments. This will be a further development of our research.

We also remark that in in our experiments we use GLCM computed with distance
1, because our aim is to show how the invariant moments can
be more discriminative if extracted from a different image representation,
rather than obtaining perfect texture classification results.
Indeed, in a recent work \cite{DiRuberto2017b} we demonstrated that a larger
distance may lead to better classification performances.
Also, GLCM computed with distance 1 are not well suited to characterize coarse
textures or even objects, both present on Vectorial and ALOT datasets (see
Figure~\ref{fig:textures}), but performs better in describing fine textures and
close patterns.

\section{Conclusions}\label{sec:conclusion}

Among all types of moments, orthogonal moments present the peculiar property of
being characterized both by small information redundancy and by high
discriminative power. For these reasons, they are widely used in different 
applications: shape-based image retrieval, edge detection, and as a feature set
in pattern recognition and in biomedical image analysis.

In this work we compute orthogonal moments by employing recurrence
relations for Legendre and discrete Chebyshev polynomials, taking into account
not only speed, but also accuracy in computation. Indeed, errors in moments
approximation may affect their discriminative ability in classification.
As it can be expected, using the described approach we achieve a great
reduction in the computational complexity and a more stable computation, with
respect to closed form representation.

We present the results of numerical experiments to assess the discriminative
power of orthogonal moments in texture analysis tasks. To this end, we use six
different databases of texture images. The Legendre polynomials and the
discrete Chebyshev polynomials computed by recursion formulas appear to be much
more effective than the Legendre polynomials computed by the closed form
representation, when the moments are extracted from the input image, and their
accuracy slowly decrease when the number of moments increases.
Additionally, the moments computed from the GLCM of the input image achieve a
high accuracy in classification, which improves as the number of moments gets
larger. Finally, we found that the use of weighted moments produces better
results in some situations; this will be the object of future work.

To further highlight the performance obtained by this approach, we compare
these results with some widely used state-of-the-art descriptors for texture
classification: rotation invariant GLCM, rotation invariant LBP, and CNN
features. Outperforming the alternative approaches in 3 datasets out of 6, and
producing equivalent results on a fourth collection of images, the orthogonal
moments reveal themselves as powerful and competitive descriptors for texture
analysis and classification.

\bibliographystyle{plain}
\bibliography{biblio}

\end{document}